\documentclass[11pt, leqno,twoside]{article}
\usepackage{amssymb}
\usepackage{amsmath}
\usepackage{amsthm}
\usepackage{graphicx}

\usepackage[active]{srcltx}

\allowdisplaybreaks

\def\rr{{\mathbb R}}
\def\rn{{{\rr}^n}}

\def\cd{{\mathcal D}}

\def\cm{{\mathcal M}}

\def\cb{{\mathcal B}}

\def\fz{\infty}
\def\az{\alpha}
\def\dist{{\mathop\mathrm{\,dist\,}}}
\def\loc{{\mathop\mathrm{\,loc\,}}}
\def\ball{{\mathop\mathrm{\,ball\,}}}
\def\lip{{\mathop\mathrm{\,Lip}}}
\def\llc{{\mathop\mathrm{\,LLC}}}
\def\ccap{{\mathop\mathrm{\,Cap}}}

\def\dz{\delta}

\def\ez{\epsilon}

\def\bz{\beta}

\def\gz{{\gamma}}

\def\boz{{\Omega}}

\def\vz{\varphi}
\def\tz{\theta}

\def\wz{\widetilde}

\def\ls{\lesssim}
\def\gs{\gtrsim}

\def\bint{{\ifinner\rlap{\bf\kern.35em--}
\int\else\rlap{\bf\kern.45em--}\int\fi}\ignorespaces}

\def\bbint{{\ifinner\rlap{\bf\kern.25em--}
\hspace{0.078cm}\int\else\rlap{\bf\kern.45em--}\int\fi}\ignorespaces}

\def\lp{{L^{p}(\boz)}}

\def\diam{{\mathop\mathrm{\,diam\,}}}

\def\r{\right}
\def\lf{\left}

\newtheorem{thm}{Theorem}[section]
\newtheorem{lem}{Lemma}[section]
\newtheorem{prop}{Proposition}[section]
\newtheorem{rem}{Remark}[section]
\newtheorem{cor}{Corollary}[section]
\newtheorem{defn}{Definition}[section]

\numberwithin{equation}{section}

\begin{document}

\arraycolsep=1pt

\title{\Large\bf  Haj\l asz-Sobolev Imbedding and Extension
\footnotetext{\hspace{-0.35cm}
\noindent{2000 {\it Mathematics Subject Classification:}} 46E35
\endgraf  {\it Key words and phases:}  Haj\l asz-Sobolev space, Haj\l asz-Sobolev extension,
Haj\l asz-Sobolev imbedding, Triebel-Lizorkin space, weak cigar domain, uniform domain, local linear connectivity
\endgraf  Yuan Zhou was supported by the Academy of Finland grant 120972.
}}
\author{Yuan Zhou}
\date{ }
\maketitle

\begin{center}
\begin{minipage}{13.5cm}\small
{\noindent{\bf Abstract}\quad
The author establishes some geometric criteria for a Haj\l asz-Sobolev $\dot M^{s,\,p}_\ball$-extension
(resp. $\dot M^{s,\,p}_\ball$-imbedding)
 domain of $\rn$ with $n\ge2$,  $s\in(0,\,1]$ and $p\in[n/s,\,\fz]$ (resp. $p\in(n/s,\,\fz]$).
In particular, the author proves that a bounded finitely connected planar domain $\boz$ is a weak $\az$-cigar domain with
$\az\in(0,\,1)$ if and only if
$\dot F^s_{p,\,\fz}(\rr^2)|_\boz=\dot M^{s,\,p}_\ball(\boz)$ for some/all $s\in[\az,\,1)$ and $p=(2-\az)/(s-\az)$, where
$\dot F^s_{p,\,\fz}(\rr^2)|_\boz$ denotes the restriction of the Triebel-Lizorkin space $\dot F^s_{p,\,\fz}(\rr^2) $
on $\boz$.
}
\end{minipage}
\end{center}

\medskip

\section{Introduction\label{s1}}

\hskip\parindent

Let $X(\boz)$ and $Y(\boz)$ be  function spaces defined on a domain  $\boz\subset \rn$.
Then $\boz$ is called an $X$-extension domain if $X(\boz)=X(\rn)|_\boz$ with equivalent norms,
where $X(\rn) |_\boz\equiv\{u|_\boz:\ u\in X(\rn)\}$ and for $v\in X(\rn)|_\boz$,
$\|v\|_{X(\rn) |_\boz}\equiv\inf\|u\|_{X(\rn)}$ with the infimum taken over all  $u\in X(\rn)$ such that $u|_\boz=v$.
Also $\boz$ is said to support an imbedding
from  $X(\boz)$ to $Y(\boz)$ if $X(\boz)$ is a subset of $ Y(\boz)$
and for all $u\in X(\boz)$,
$\|u\|_{Y(\boz)}\le C \|u\|_{X(\boz)}$ with constant $C$ independent of $u$.
Moreover, we always denote by $\dot W^{1,\,p}(\boz)$
(resp. $W^{1,\,p}(\boz)$) with $p\in(1,\,\fz]$ the homogeneous (resp. inhomogeneous) Sobolev space.
The other notions for domains, such as uniform domain,  weak $\az$-cigar domain, regular domain,
the LLC property and the slice property, will be explained in Section \ref{s2}.

For the history of  geometric properties of Sobolev extension and imbedding domains
see, for example, \cite{gv,j81,m85,gm85,gm85b,k90,bk95,bk96,k98,bs,s09,s10,hkt08,hkt08b} and their references.
In particular, it was proved by Jones \cite{j81} that
a bounded simply connected domain $\boz\subset\rr^2$ is a uniform domain
if and only if it is a $W^{1,\,2}$-extension domain;
and by Gehring and Martio \cite{gm85b} that the
$ W^{1,\,n}$-extension domain has the LLC property;
see also \cite{k90,g82,gr83,gv,vgl}.
On the other hand, Gehring and Martio \cite{gm85} proved that
for $\az\in(0,\,1]$, $\boz$ is a weak $\az$-cigar domain
 if and only if it is a local $\lip_\az$-extension  domain.
Moreover, let $p\in(n,\,\fz)$ and $\az=(p-n)/(p-1)$.
It was proved by Buckley and Koskela  \cite{bk96} that a weak $\az$-cigar domain
always supports an imbedding from $\dot W^{1,\,q}(\boz)$ into $\dot C^{1-n/q}(\overline \boz)$   for all $q\in[p,\,\fz)$;
and by Koskela \cite{k98} that a weak $\az$-cigar domain is a $ W^{1,\,q}$-extension domain for all $q\in(p,\,\fz)$,
which was further improved by Shvartsman \cite{s10} to
all $q\in(p^\ast,\,\fz)$ with some $p^\ast\in(n,\,p)$.
Conversely,
with the additional assumption that $\boz$ has the slice property,
Buckley and Koskela \cite{bk96} proved that
if $\boz$ supports an imbedding from  $\dot W^{1,\,p}(\boz)$ into $\dot C^{1-n/p}(\overline \boz)$,
then it is a weak $\az$-cigar domain.
We also refer the reader to \cite{bs} for more criteria for $\dot W^{1,\,p}$-imbedding domains,
where they reduce the  slice property to some weak slice properties.

It was noticed by Haj\l asz \cite{h96} that the simple pointwise inequality
\begin{equation}\label{e1.1}
|u(x)-u(y)|\le |x-y|^s[g(x)+g(y)]
\end{equation}
can be used to characterize Sobolev functions $u$ when $s=1$. More generally, for $s\in(0,\,1]$
and measurable function $u$,
denote by $\cd^s(u)$ the collection of all nonnegative measurable functions $g$
such that \eqref{e1.1} holds
for all $x,\ y\in\boz\setminus E$, where
$E\subset\boz$ satisfies $|E|=0$.
We also denote by $\cd_\ball^s(u)$ the collection of
all nonnegative measurable functions $g$ such that
 \eqref{e1.1} holds for all $x,\ y\in\boz\setminus E$
satisfying $|x-y|<\frac12\dist(x,\,\partial\boz)$.

\begin{defn}\rm\label{d1.1}
Let $s\in (0,1]$ and $p\in(0,\,\fz]$.

(i) The homogeneous Haj\l asz  space
 $\dot M^{s,\,p}(\boz)$ is the
space of all measurable functions $u$ such that
$$\|u\|_{\dot M^{s,\,p} (\boz)}\equiv\inf_{g\in \cd^s (u)}\|g\|_\lp<\fz.$$

(ii) The Sobolev-type Haj\l asz space $\dot M^{s,\,p}_\ball(\boz)$  is the
space of all measurable functions $u$ such that
$$\|u\|_{\dot M^{s,\,p}_\ball(\boz)}\equiv\inf_{g\in \cd^s_\ball(u)}\|g\|_\lp<\fz.$$

Moreover, we set
$M^{s,\,p}(\boz)\equiv L^p(\boz)\cap \dot M^{s,\,p}(\boz)$
with $\|u\|_{M^{s,\,p}(\boz)}\equiv\|u\|_{\dot M^{s,\,p}(\boz)}+\|u\|_{L^p(\boz)}$
for all $u\in M^{s,\,p}(\boz)$, and similarly define $M^{s,\,p}_\ball(\boz)$.
\end{defn}

Obviously, for all $s\in(0,\,1]$ and $p\in(0,\,\fz]$,
$\dot M^{s,\,p}(\boz)\subset \dot M^{s,\,p}_\ball(\boz)$.
Conversely, if $\boz$ is a uniform domain,
then $\dot M^{s,\,p}_\ball(\boz)=\dot M^{s,\,p}(\boz)$ for all $s\in(0,\,1]$ and $p\in(n/(n+s),\,\fz]$;
 see \cite[Theorem 19]{ks08} and also \cite[Lemma 14]{hm97}. But, generally, we cannot expect that $\dot M^{s,\,p}(\boz)=\dot M^{s,\,p}_\ball(\boz)$.
For example,  this fails when $\boz=B(0,\,1)\setminus \{(x,\,0):\,x\ge0\}\subset\rr^2.$

Haj\l asz-Sobolev spaces are closely related to the classical (Hardy-)Sobolev and Triebel-Lizorkin spaces.
In fact, it was proved in \cite{h96,ks08} that
$\dot W^{1,\,p}(\boz)=\dot M^{1,\,p}_\ball(\boz)$ for $p\in(1,\,\fz]$
 and  $\dot H^{1,\,p}(\boz)=\dot M^{1,\,p}_\ball(\boz)$ for $p\in(n/(n+1),\,1]$,
 which together with \cite{t83} implies that
$\dot M^{1,\,p}(\rn)=\dot M^{1,\,p}_\ball(\rn)=\dot F^1_{p,\,2}(\rn)$ for all $p\in(n/(n+1),\,\fz]$,
while for all $s\in(0,\,1)$ and $p\in(n/(n+s),\,\fz]$,
 $\dot M^{s,\,p}(\rn)=\dot M^{s,\,p}_\ball(\rn)=\dot F^s_{p,\,\fz}(\rn)$
as proved in \cite{y03,kyz09,kyz09c}.
Here and in what follows,
we always denote
by $\dot H^{1,\,p}(\boz)$ with $p\in(0,\,1]$ the Hardy-Sobolev space as in \cite{m90},
and by  $\dot F^s_{p,\,q}(\rn)$ with $s\in\rr$ and $p,\,q\in(0,\,\fz]$
the homogeneous Triebel-Lizorkin spaces as in \cite{t83}.

Recently, it was proved in \cite{hkt08b} (see  \cite{hkt08,s06} and also Lemma \ref{l4.1} below) that for $p\in(1,\,\fz)$,
$\boz$ is a $W^{1,\,p}$-extension if and only if $\boz$ is regular (see Definition \ref{d2.5})
and $\dot W^{1,\,p}(\boz)=\dot M^{1,\,p}(\boz)$ (namely, $\dot M_\ball^{1,\,p}(\boz)=\dot M^{1,\,p}(\boz)$),
while $\boz$ is regular if and only if
 $\boz$ is an $M^{1,\,p}$-extension domain.
Some characterizations of the restriction of Besov and Triebel-Lizorkin spaces on  regular domains
were also established by Shvartsman \cite{s07}.
Recall that it is an interesting subject to establish
some intrinsic characterizations of $\dot F^s_{p,\,q}(\rn)|_\boz$,
the restriction of the Triebel-Lizorkin space  $\dot F^s_{p,\,q}(\rn)$ on the domain $\boz$;
see \cite {r99,r00,t83,t02} for more discussions.
In particular, some intrinsic characterizations of the restriction of
Triebel-Lizorkin spaces on Lipschitz domains were established
by Rychkov \cite{r99,r00} and Triebel \cite{t02}.

In what follows, $\boz$ is called an $\dot M^{s,\,p}_\ball$-imbedding domain
if it supports an imbedding from $\dot M^{s,\,p}_\ball(\boz)$ to
$\dot M^{s-n/p,\,\fz} (\boz)$  with $s\in(0,\,1]$ and $p\in(n/s,\,\fz]$.
We define $ M^{s,\,p}_\ball$-imbedding domains similarly.

The main purpose of this paper is to establish some geometric criteria for $\dot M^{s,\,p}_\ball$-extension
(resp. $\dot M^{s,\,p}_\ball$-imbedding)
 domains of $\rn$ with $n\ge2$,  $s\in(0,\,1]$ and $p\in[n/s,\,\fz]$ (resp. $p\in(n/s,\,\fz]$).
In particular, we prove  that a bounded simply connected planar domain $\boz$ is a weak $\az$-cigar domain with
$\az\in(0,\,1)$ if and only if
$\dot F^s_{p,\,\fz}(\rr^2)|_\boz=\dot M^{s,\,p}_\ball(\boz)$ for some/all $s\in[\az,\,1)$ and $p=(2-\az)/(s-\az)$.

More precisely, we first obtain the following conclusion
by using  some ideas from \cite{gm85,k90,hek98}  and introducing a
capacity associated to  $\dot M^{s,\,n/s}_\ball(\boz)$.
See  Section \ref{s3} for its proof.

\begin{thm}\label{t1.1}
If $\boz$ is a bounded $\dot M^{s,\,n/s}_\ball$-extension domain for some $s\in(0,\,1]$,
then $\boz$ has the LLC property.
\end{thm}

Recall that if a bounded simply connected planar domain, or a bounded
domain of $\rn$ with $n\ge2$ that is quasiconformally equivalent to a uniform domain,
has the LLC property, then it is a uniform domain; see \cite{k90}.
We also recall that $\dot F^{s}_{p,\,\fz}(\rn)=\dot M^{s,\,p}_\ball(\rn)$ for all $s\in(0,\,1)$
and $p\in(n/(n+s),\,\fz]$;
see \cite{y03} and also \cite{kyz09}.
 Then as a corollary to  Theorem \ref{t1.1},
we have the following conclusion.

\begin{cor}\label{c1.1}
Let $\boz$ be a bounded simply connected planar domain, or a bounded
domain of $\rn$ with $n\ge2$ that is quasiconformally equivalent to a uniform domain.
Then the following are equivalent:

(i) $\boz$ is a uniform domain;

(ii)  $\boz$ is an $\dot M^{s,\,n/s}_\ball$-extension domain  for some/all $s\in(0,\,1]$;

(iii)  $\dot F^{s}_{n/s,\,\fz}(\rn)|_\boz=\dot M^{s,\,n/s}_\ball(\boz)$
for some/all $s\in(0,\,1)$.
\end{cor}

When $p\in(n/s,\,\fz)$, we also establish the following geometric characterizations,
which generalizes \cite[Theorem 4.1]{bk96} and \cite[Theorem 1.1]{s10} to Haj\l asz-Sobolev spaces.
See Section \ref{s4} for its proof, which uses some ideas from \cite{bk96,k98,s10}, in particular,
uses Theorems \ref{t4.1} and \ref{t4.2} below, and the weak self-improving property
of a weak cigar domain established by Shvartsman in \cite[Theorem 1.4]{s10}
(see also Proposition \ref{p4.1} below).

\begin{thm}\label{t1.2}
(i) Let $\az\in(0,\,1)$ and $\boz\subset\rn$ be a bounded weak $\az$-cigar domain.
Then for all $s\in(\az,\,1]$ and $p\in[(n-\az)/(s-\az),\,\fz)$,
 $\boz$ is an $\dot M^{s,\,p}_\ball$-extension domain and, especially, an
 $\dot M^{s,\,p}_\ball $-imbedding domain.

(ii)
Let $s\in(0,\,1]$, $p\in(n/s,\,\fz)$ and $\az\in[(ps-n)/(p-1),\,1]$.
If $\boz$ is an bounded $\dot M^{s,\,p}_\ball$-extension or $\dot M^{s,\,p}_\ball$-imbedding domain
having the slice property,
then $\boz$ is a weak $\az$-cigar domain.
\end{thm}

At the endpoint case $p=\fz$, as proved by Gehring and Martio \cite{gm85}, a bounded domain
$\boz$ is a weak $\az$-cigar domain with $\az\in(0,\,1]$ if and only if it is an $\dot M^{\az,\,\fz}_\ball$-extension domain,
and if and only if it is an $\dot M^{\az,\,\fz}_\ball$-imbedding domain,
where $\dot M^{\az,\,\fz}(\boz)$ and $\dot M^{\az,\,\fz}_\ball(\boz)$
 coincide with  $\lip_\az(\boz)$ and  $\loc\lip_\az(\boz)$  as in \cite{gm85},
respectively.
Recall that a bounded simply connected planar domain, or a bounded
domain of $\rn$ with $n\ge2$ that is quasiconformally equivalent to a uniform domain,
always has the slice property (see \cite{bk96}).
Then, as a corollary to Theorem \ref{t1.2} and \cite{gm85}, we have the following conclusion,
which together with Corollary \ref{c1.1} gives an intrinsic characterization
of the restriction of the Triebel-Lizorkin space
$\dot F^s_{p,\,\fz}(\rn)|_\boz$ for a class of domains $\boz$.

\begin{cor}\label{c1.2}
Let $\az\in(0,\,1)$ and $\boz$ be a bounded simply connected planar domain, or a  bounded
 domain of $\rn$ with $n\ge2$ that is quasiconformally equivalent to a uniform domain.
Then the following are equivalent:

(i) $\boz$ is a weak $\az$-cigar domain;

(ii)  $\dot F^{s}_{p,\,\fz}(\rn)|_\boz=\dot M^{s,\,p}_\ball(\boz)$
for some/all $s\in[\az,\,1)$ and $p=(n-\az)/(s-\az)$;

 (iii) $\boz$ is an $\dot M^{s,\,p}_\ball $-extension domain
for some/all $s\in[\az,\,1]$ and $p=(n-\az)/(s-\az)$;

(iv) $\boz$ is  an $\dot M^{s,\,p}_\ball $-imbedding  domain
for some/all $s\in[\az,\,1]$ and $p=(n-\az)/(s-\az)$.
\end{cor}

Finally, let $\az\in(0,\,1)$ and $\boz$ be a bounded weak $\az$-cigar domain,
namely, bounded $\az$-subhyperbolic domain as in \cite{s10}.
Then, with the aid of its weak self-improving property established in \cite[Theorem 1.5]{s10}
(see also Proposition \ref{p4.1} below),
Shvartsman \cite[Theorem 1.1]{s10} proved that $\boz$ is a $W^{1,\,p}$-extension domain
for every $p\in((n-\az^\ast)/(1-\az^\ast),\,\fz)$ with
$\az^\ast\in(0,\,\az)$ as in Proposition \ref{p4.1}.
Following this and \cite[Theorems 1.1]{s10} with taking
$p\in((n-\az^\ast)/(1-\az^\ast),\,(n-\az)/(1-\az))$, 
if $\boz$ is also a finitely connected planar domain,
then Shvartsman \cite[p.\,2210]{s10} pointed out 
that $\boz$ is a weak $\tau$-cigar domain with $\tau\in(\az^\ast,\,\az)$,
while such a self-improving property for an arbitrary weak cigar domain is still unknown.
However, as Shvartsman informed me, when $p\in((n-\az^\ast)/(1-\az^\ast),\,(n-\az)/(1-\az))$,
the proof of \cite[Theorem 1.1]{s10} actually proved 
a weak form of \cite[Theorem 1.1]{s10}: there exists a continuous linear extension operator from 
$ W^{1,\,p}(\boz)\cap W^{1,\, (n-\az)/(1-\az)}(\boz)$ to $ W^{1,\,p}(\boz)$. 
But, the following conclusion shows that $\boz$ does have
  the self-improving property when it is also a finitely connected planar domain.
The point is that, as observed in Remark \ref{r4.3},
if $\boz$ is a finitely connected planar domain or more generally, $\boz$ has the slice property,
then the $\dot M^{s,\,p}_\ball(\boz)$-imbedding
required in Theorem \ref{t1.2}(ii) can be reduced to a weaker one,
which is already obtained in the proof of Theorem \ref{t1.2}(i)
with the aid of Shvartsman \cite[Theorem 1.5]{s10} or 
which, when $s=1$, can also be deduced from the above weak form of \cite[Theorem 1.1]{s10} 
as Shvartsman informed me.
See Section \ref{s4} for more details.

\begin{thm}\label{t1.3}
Let $n\ge 2$ and  $\boz\subset\rn$ be a  finitely connected bounded planar domain or
a bounded domain that is quasiconformally equivalent to a uniform domain, or more generally,   
 be a bounded domain satisfying the slice property.
If  $\boz $ is a  weak $\az$-cigar domain  with some $\az\in(0,\,1)$,
then  $\boz$ is a weak $\az^\ast$-cigar domain with $\az^\ast\in(0,\,\az)$,
and hence for all $s\in(\az^\ast,\,1]$ and $p\in[(n-\az^\ast)/(s-\az^\ast),\,\fz)$,
$\boz$ is an $\dot M^{s,\,p}_\ball$-extension domain and, especially, an
$\dot M^{s,\,p}_\ball $-imbedding domain.
\end{thm}

\section{Preliminaries}\label{s2}

In this section, we recall some  notions and basic properties of
domains and  Haj\l asz spaces.
We begin with the notion of a uniform domain.

\begin{defn}\label{d2.1}\rm
A domain
 $\boz\subset\rn$ is called a uniform domain if there exists a
positive constant $  C$ such that for all $x,\,y\in\boz$, there exists a rectifiable curve $\gz:\,[0,\,T]\to\boz$,
parameterized by the arc length, with $\gz(0)=x$ and $\gz(T)=y$, and satisfying  that $T\le   C|x-y|$ and
\begin{equation}\label{e2.1}
\bigcup_{t\in[0,\,T]}B\lf(\gz(t),\,\frac1{  C}\min\{t,\,T-t\} \r)\subset\boz.
\end{equation}
\end{defn}

Closely related to the concept of a
uniform domain, Gehring \cite{g77} introduced the notion of  linear local connectivity.
\begin{defn}\label{d2.2}\rm
A domain $\boz\subset\rn$  is said to have the  linearly locally connectivity (for short, LLC) property if
there exists a constant $b \in(0,\,1]$ such that
 for all $z\in\rn$ and $r>0$,

\noindent$\llc(1)$ \quad points in $\boz\cap B(z,\,r)$ can be joined in $\boz\cap B(z,\,r/b)$;

\noindent$\llc(2)$ \quad points in $\boz\setminus B(z,\,r)$ can be joined in $\boz\setminus B(z,\,br)$.
\end{defn}

It is known that each uniform domain has the LLC property.
Conversely, assume that $\boz$ is
a bounded simply connected planar domain, or
 a bounded domain of $\rr^n$ with $n \ge 3$ that is quasiconformally equivalent to a uniform domain.
If $\boz$ has the LLC property, then it is a uniform domain; see \cite{v84} and also \cite{k90}.

Now we recall the notion of a weak  cigar domain; see \cite{gm85,bk96}.

\begin{defn}\label{d2.3}\rm  Let  $\az\in(0,\,1]$.
Then a domain $\boz\subset\rn$ is called a weak $\az$-cigar domain
if there exists a
positive constant $  C$ such  that for every pair of points $x,\,y\in\boz$,
there exists a rectifiable curve $\gz\subset \boz$ joining $x$ and $y$, and satisfying
$$
 \int_\gz [d(z,\,\boz^\complement)]^{\az-1}\,|dz|\le C|x-y|^\az.
$$
\end{defn}

Notice that the class of weak $1$-cigar domains coincides with the class of quasiconvex domains,
and a bounded weak $\az$-cigar domain is a weak $\bz$-cigar domain for all $\bz\in(\az,\,1]$;
see \cite{bk96} for  details. Moreover, a uniform  domain is also a weak $\az$-cigar domain for all $\az\in(0,\,1]$.
.

The following slice property was introduced by Buckley and Koskela \cite{bk96}.
In what follows, for every rectifiable curve $\gz$, we denote its length by $\ell(\gz)$.

\begin{defn}\label{d2.4}\rm
A domain $\boz$ has a slice property with respect to   $C>1$
if for every pair of points $x,\,y \in\boz$, there exists  a rectifiable curve $\gz:[0,\,1]\to\boz$ with
$\gz(0)=x $ and  $\gz(1)=y$, and pairwise disjoint collection of open subsets
$\{S_i\}_{i=0}^j$, $j\ge0$,  of $\boz$ such that

(i) $x\in S_0$, $y\in S_j$ and $x $ and $y$ are in different components of $\boz\setminus \overline{S_i}$ for $0<i<j$;

(ii) if $F\subset\subset  \boz$ is a curve containing both $x$ and $y$, and $0<i<j$, then $\diam(S_i)\le C\ell(F\cap S_i)$;

(iii) for $0\le t\le 1$, $B(\gz(t),\,C^{-1}d(\gz(t),\,\boz^\complement))\subset\cup_{i=0}^jS_i$;

(iv) if $0\le i\le j$, then $\diam S_i\le C d(z,\,\boz^\complement)$ for all $z\in\gz_i\equiv\gz\cap S_i$;
also, there exists $x_i\in S_i$ such that $x_0=x$, $x_j=y$ and
$B(x_i,\,C^{-1}d(x_i,\,\boz^\complement))\subset S_i$.
\end{defn}

It was proved by Buckley and Koskela \cite{bk96} that every simply connected domain in $\rr^2$, or
every domain in $\rr^n$ with $n \ge 3$ that is quasiconformally equivalent to a uniform domain,
has the  slice property  as in Definition \ref{d2.4}.

We also recall the notion of the regularity of a domain.
\begin{defn}\label{d2.5}
A domain $\boz\subset\rn$ is  regular if there exist  positive constants $\tz$ and $C$ such that
for all $x\in\overline\boz$ and $r\in(0,\,\tz)$,
$|B(x,\,r)\cap\boz|\ge C|B(x,\,r)|.$
\end{defn}

We point out that the regularity of $\boz$ does not depend on the choice of $\tz$ and $C$ in the following sense:
if   $\boz $ is  regular with   $\tz$ and $C$, then for any $\tz'\in(0,\,\fz)$, there exists a constant $C'$ such that
 $\boz $ is  regular with   $\tz'$ and $C'$.

The following lemma established in \cite{ks08}
will be useful in the following proofs.
In what follows, for every $\rho\in(0,\,\fz]$, similarly to $\cd_\ball^s(u)$,
we denote by $\cd_\ball^{s,\,\rho}(u)$ the collection of
all measurable functions $g$ such that
 \eqref{e1.1} holds for all $x,\ y\in\boz\setminus E$
satisfying $|x-y|<\rho\dist(x,\,\partial\boz)$.
Notice that $\cd_\ball^s(u)= \cd_\ball^{s,\,1/2}(u)$ and
$\cd^s(u)= \cd_\ball^{s,\,\fz}(u)$.

\begin{lem}\label{l2.1}
Let $s\in(0,\,1]$ and $p\in(n/(n+s),\,\fz]$. Then $u\in\dot M^{s,\,p}_\ball(\boz)$ if and only if
there exists a $\rho\in(0,\,1)$ such that
$\inf_{g\in\cd^{s,\,\rho}_\ball(u)}\|g\|_{L^p(\boz)}<\fz$.
Moreover, for given $\rho$,  there exists a positive constant $C$ such that
for all $u\in\dot M^{s,\,p}_\ball(\boz)$,
$$C^{-1}\|u\|_{\dot M^{s,\,p}_\ball(\boz)}\le \inf_{g\in\cd^{s,\,\rho}_\ball(u)}\|g\|_{L^p(\boz)}\le C\|u\|_{\dot M^{s,\,p}_\ball(\boz)}.$$
\end{lem}

Finally, we state some conventions. Throughout the paper,
we denote by $C$ a positive constant which is independent
of the main parameters, but which may vary from line to line.
Constants with subscripts, such as $C_0$, do not change
in different occurrences. The symbol $A\ls B$ or $B\gs A$
means that $A\le CB$. If $A\ls B$ and $B\ls A$, we then
write $A\sim B$.
For any locally integrable function $f$,
we denote by $\bbint_E f\,d\mu$ the average
of $f$ on $E$, namely, $\bbint_E f\,d\mu\equiv\frac 1{|E|}\int_E f\,dx$.

\section{Proof of Theorem \ref{t1.1} }\label{s3}

To prove Theorem \ref{t1.1}, we need the following capacity.

\begin{defn}\label{d3.1}
Let $s\in(0,\,1]$. For every pair $E,\,F\subset\boz$ of disjoint continua, define the capacity associated to $\dot M^{s,\,n/s}_\ball(\boz)$ by
$$\ccap_{\dot M^{s,\,n/s}_\ball}(E,\,F,\,\boz)\equiv\inf_{u\in\Delta_s(E,\,F,\,\boz)}\|u\|^{n/s}_{\dot M^{s,\,n/s}_\ball(\boz)},$$
where $\Delta_s(E,\,F,\,\boz)$ denotes the collection of all continuous functions $u\in\dot M^{s,\,n/s}_\ball(\boz)$
with $u(x)=0$ for all $x\in E$ and $u(x)=1$ for all $x\in F$.
\end{defn}

Obviously,  for every pair $E,\,F\subset\boz$  of disjoint continua
and every pair $\wz E,\,\wz F\subset\rn$ of disjoint continua satisfying
 $E\subset \wz E$ and $F\subset \wz F$, we have
\begin{equation}\label{e3.1}
\ccap_{\dot M^{s,\,n/s}_\ball}(E,\,F,\,\boz)\le \ccap_{\dot M^{s,\,n/s}_\ball}(\wz E,\,\wz F,\,\rn).
\end{equation}
Moreover, a reverse inequality also follows for  $\dot M^{s,\,n/s}_\ball$-extension domains
by modifying the proof of \cite[Theorem 2.2]{k90}.  We omit the details.

\begin{lem}\label{l3.1}
If $\boz$ is an $\dot M^{s,\,n/s}_\ball$-extension domain, then
there exists a positive constant $C$ such that
for every pair $E,\,F\subset\boz$ of disjoint continua,
$$\ccap_{\dot M^{s,\,n/s}_\ball}(E,\,F,\,\rn)\le C\ccap_{\dot M^{s,\,n/s}_\ball}(E,\,F,\,\boz).$$
\end{lem}

The following property of the capacity plays an important role in the proof of Theorem \ref{t1.1}
and it is proved by using some ideas of \cite[Theorem 5.9]{hek98}.

\begin{lem}\label{l3.2}
 Let $s\in(0,\,1]$ and $\dz\in(0,\,\fz)$. There exists a positive constant $C$ such that
for every pair $E,\,F\subset\rn$ of disjoint continua, if
$\min\{\diam E,\,\diam F\}\ge \dz \dist(E,\,F),$
then $\ccap_{\dot M^{s,\,n/s}_\ball}(E,\,F,\,\rn)\ge C.$
\end{lem}

\begin{proof}
Notice that if continua $\wz F\subset F$
and $\wz E\subset E$, then $ \Delta_s(E,\,F,\,\rn)\subset  \Delta_s(\wz E,\,\wz F,\,\rn)$
and thus $$\ccap_{\dot M^{s,\,n/s}_\ball}(E,\,F,\,\rn)\ge \ccap_{\dot M^{s,\,n/s}_\ball}(\wz E,\,\wz F,\,\rn).$$
So without loss of generality, we may assume that $\diam E=\diam F\ge \dz \dist(E,\,F)$.

Fix $x_0\in F$ and $r\equiv (2+\dz)\diam E$. Then $E,\,F\subset B(x_0,\,r)$. Let $u\in  \Delta_s (E,\,F,\, \rn)$
and, without loss of generality, assume that $u_{B(x_0,\,r)}\le1/2$.
Then for every $x\in F$ and $g\in\cd_\ball^s(u)\cap L^{n/s}(\boz)$, we have
\begin{eqnarray*}
 \frac12\le|u(x)-u_{B(x_0,\,r)}|&&\le \sum_{i=-1}^\fz|u_{B(x,\,2^{-i}r)}-u_{B(x,\,2^{-i-1}r)}|+|u_{B(x,\,2r)}-u_{B(x_0,\,r)}|\\
&&\le \sum_{i=-1}^\fz (2^{-i}r)^s \lf(\bint_{B(x,\,2^{-i}r)}[g(z)]^{n/s}\,dz\r)^{s/n}\\
&&\le \sum_{i=-1}^\fz (2^{-i}r)^{s/n}\lf(\frac1{2^{-i}r}\int_{B(x,\,2^{-i}r)}[g(z)]^{n/s}\,dz\r)^{s/n}\\
&&\ls  \sup_{0<t\le2r}\lf(\frac rt\int_{B(x,\,t)}[g(z)]^{n/s}\,dz\r)^{s/n},
\end{eqnarray*}
which implies that there exists $t_x\in(0,\,2r]$ such that
$$t_x\ls r\int_{B(x,\,t)}[g(z)]^{n/s}\,dz.$$
By the Vitali covering lemma, we can find points $\{x_i\}_i\subset F$ such that
$\{B(x_i,\,t_{x_i})\}_i$ are pairwise disjoint and $F\subset \cup_i 5B(x_i,\,t_{x_i})$.
Thus,
$$\diam F\le\sum_i 10t_{x_i}\ls r\sum_i\int_{B(x_i,\,t_{x_i})} [g(z)]^{n/s}\,dz\ls r\int_{B(x_0,\,3r)}[g(z)]^{n/s}\,dz,$$
which yields that $\|u\|_{\dot M^{s,\,p}_\ball(\boz)}\gs 1$ and thus finishes the proof of Lemma \ref{l3.2}.
\end{proof}

\begin{proof}[Proof of Theorem \ref{t1.1}.]
We first prove that $\boz$ has the LLC(2) property.
Let $x_1,\,x_2\in B(x_0,\,r)\cap\boz$ for some $x_0\in\rn$ and $r>0$.
Suppose that $x_1$ and $x_2$ are not in the same component of $\boz\setminus B(x_0,\,b_0r)$
with $b_0\in(0,\,1/4)$. It then suffices to prove that $b_0$ is bounded from below.
To this send, we choose a rectifiable curve $\gz\subset \boz$ joining $x_1$ and $x_2$,
and denote by $F_i$ the component of $\gz\cap(\boz\setminus B(x_0,\,r/2))$ containing $x_i$ for $i=1,\,2$.
Obviously, $\diam F_i\ge r/2\ge\dist(F_1,\,F_2)/4$ for $i=1,\,2$. Then by \eqref{e3.1} and Lemmas \ref{l3.1} and \ref{l3.2},
we have
\begin{equation}\label{e3.2}
 \ccap_{\dot M^{s,\,n/s}_\ball}(F_1,\,F_2,\,\boz)\sim  \ccap_{\dot M^{s,\,n/s}_\ball}(F_1,\,F_2,\,\rn)\gs1.
\end{equation}
To estimate $ \ccap_{\dot M^{s,\,n/s}_\ball}(F_1,\,F_2,\,\boz)$ from above, for all $x\in\rn$, define
$$u(x)\equiv\lf\{\begin{array}{ll}
1,& x\in \boz \cap B(x_0,\,b_0r);\\
 \lf(\log\frac{1}{2b_0}\r)^{-1} \lf(\log\frac{r}{2|x-x_0|}\r), \quad&      x\in \boz\cap (B(x_0,\,r/2)\setminus B(x_0,\,b_0r));\\
  0,& x\in\boz\setminus B(x_0,\,r/2)
                 \end{array}
 \r.$$
and
$$g(x)\equiv \frac1{|x-x_0|^{s}}\lf(\log\frac{1}{2b_0}\r)^{-1} \chi_{\boz\cap\overline{B(x_0,\,r/2)\setminus B(x_0,\,b_0r)}}(x).$$
Then we claim that there exists a positive constant independent of $u,\,x_0,\,b_0,\,r$ such that
$Cg$ is  an element of $\cd^{s,\,1/32}_\ball(u)$.
Assume that the claim holds for the moment.  Then by Lemma \ref{l2.1},  $u\in\dot M^{s,\,p}_\ball(\boz)$ and
\begin{eqnarray*}
 \|u\|_{\dot M^{s,\,n/s}_\ball(\boz)}&&\le\|g\|_{L^{n/s}(\boz)}\\
&&\ls
\lf(\log\frac{1}{2b_0}\r)^{-1}\lf\{\int_{B(x_0,\,r/2)\setminus B(x_0,\,b_0r)} |z-x_0|^{-n}\,dz\r\}^{s/n}\\
&&\ls \lf(\log\frac{1}{2b_0}\r)^{s/n-1}.
\end{eqnarray*}
Moreover, observe that $u\in\Delta_s(F_1,\,F_2,\,\boz)$. So we have
\begin{equation}\label{e3.3}
    \ccap_{\dot M^{s,\,n/s}_\ball}(E,\,F,\,\boz)\ls
 \lf(\log\frac{1}{2b_0}\r)^{1-n/s},
     \end{equation}
which together with \eqref{e3.2} implies that $b_0\gs1$ and hence reduces the LLC(2) property of $\boz$ to proving the above claim.

To prove the above claim,  it suffices to prove that
for all $x,\,y\in\boz$ satisfying $|x-y|\le \min\{d(x,\,\boz^\complement),\,d(x,\,\boz^\complement)\}/16$,
$$|u(x)-u(y)|\ls |x-y|^s[g(x)+g(y)].$$
Fix such $x,\,y\in\boz$.
Without loss of generality, we assume that $0\le u(y)<u(x)\le1$.
Then $x\in\boz\cap \overline{B(x_0,\,r/2)}$ and by $\overline{B(x_0,\,b_0 r)}\cap\boz^\complement\ne\emptyset$,
we have $d(x,\,\boz^\complement)\le |x-x_0|+b_0r< 3r/4$.
It will not happen that $x\in \boz\cap B(x_0,\,b_0r)$ and $y\in \boz\setminus B(x_0,\, r/2)$ since, in this case,
$|x-y|\ge r/4\ge d(x,\,\boz^\complement)/3$.
If $x,\,y\in  \boz\cap \overline{B(x_0,\,r/2)\setminus B(x_0,\,b_0r)}$,
then
$$|u(x)-u(y)|=\lf(\log\frac{1}{2b_0}\r)^{-1}\lf(\log\frac{ |y-x_0|}{ |x-x_0|}\r)\ls \lf(\log\frac{r}{2b_0}\r)^{-1}\frac{|x-y|^s}{|x-x_0|^s}.$$
If $x\in \boz\cap \overline{B(x_0,\,b_0r)}$ and $y\in  \boz\cap \overline{B(x_0,\,r/2)\setminus B(x_0,\,b_0r)}$,
\begin{eqnarray*}
 |u(x)-u(y)|&&= \lf|1-\lf(\log\frac{r}{2b_0}\r)^{-1}\lf(\log\frac r{2 |y-x_0|}\r)\r|\\
&&\ls \lf(\log\frac{1}{2b_0}\r)^{-1}\frac{|y-x_0|-b_0r  }{|y-x_0| }\ls \lf(\log\frac{r}{2b_0}\r)^{-1}\frac{|x-y|^s}{|y-x_0|^s}.
\end{eqnarray*}
If $x\in  \boz\cap \overline{B(x_0,\,r/2)\setminus B(x_0,\,b_0r)}$ and $y\in \boz\setminus B(x_0,\, r/2)$,
\begin{eqnarray*}
 |u(x)-u(y)|&&= \lf(\log\frac{r}{2b_0}\r)^{-1}\lf(\log\frac r{2 |x-x_0|}\r)\\
&&\ls \lf(\log\frac{1}{2b_0}\r)^{-1}\frac{ r/2-|x-x_0|  }{|x-x_0| }\ls \lf(\log\frac{r}{2b_0}\r)^{-1}\frac{|x-y|^s}{|x-x_0|^s}.
\end{eqnarray*}
This shows the above claim and thus proves that $\boz$ has the LLC(2) property.

To prove that $\boz$ has LLC(1) property, it suffices to prove that $\boz$ is quasiconvex, namely, for every
pair $x_1,\,x_2$ of points in $\boz$, there exists a curve $\gz\subset\boz$ joining them with $\ell(\gz)\le C|x-y|$,
where the constant $C$ is independent of $x,\,y $ and $\gz$.

To this end, let $x_1,\,x_2$ be a pair of points in $\boz$.
If $|x_1-x_2|<\max\{d(x_1,\boz^\complement),\,d(x_2,\boz^\complement)\}$,
then the line segment joining $x_1$ and $x_2$ is the desired curve.
Assume that $|x_1-x_2|\ge \max\{d(x_1,\boz^\complement),\,d(x_2,\boz^\complement)\}$.
Let $\gz^{(0)}\subset\boz$ be a curve joining $x_1$ and $x_2$,
and let $F_i$ be the component containing $x_i$ of $\gz^{(0)}\cap B(x_i,\,
|x_1-x_2|/4)$ for $i=1,\,2$.
Notice that \eqref{e3.2} still holds by the same argument.
Moreover, there exists a positive constant $N_0>1$ independent of $x_1,\,x_2,\,\gz^{(0)},\,F_1,\,F_2$ such that
$F_1,\,F_2$ are in the same component of $\boz\cap \overline{B(x_1,\,N_0|x_1-x_2|)}$.
To see this, assume that $F_1,\,F_2$ are not in the same component of $\boz\cap \overline{B(x_1,\,N|x_1-x_2|)}$ for some $N>2$.
Then $\overline{B(x_1,\,N|x_1-x_2|)}\cap\boz^\complement\ne\emptyset$ and hence
by an argument similar to the  proof of \eqref{e3.3},  we have
$$
\ccap_{\dot M^{s,\,n/s}_\ball}(E,\,F,\,\boz)
\ls \lf(\log N\r)^{1-n/s},
$$
which means $N\ls1$ and hence shows the existence of $N_0$.
Therefore, letting $L$ be the infimum of the length of all curves joining
$F_1$ and $F_2$, we have $L<\fz$ and then define  the function
$$v(x)\equiv L^{-1}\inf_{\gz} \ell(\gz\cap B(x_1,\,N_0|x_1-x_2|))$$ for all $x\in\boz$,
where the infimum is taken over all the rectifiable curves $\gz\subset\boz$ joining $x$ and $F_1$.
For all $x\in\boz$, define $$h(x)\equiv CL^{-1}(N_0|x_1-x_2|)^{1-s}\chi_{\boz\cap\overline{B(x_1,\,N_0|x_1-x_2|)}}(x).$$
Then we claim that
there exists a positive constant $C$ independent of $v,\, L,\,x_1,\,x_2,\,N_0$ such that $Ch$ is an element of
$\cd^{s,\,1/16}_\ball(v)$. Assume that this claim holds for the moment.
Set $\wz v\equiv\min\{v,\,1\}$. Then $h$ is also a constant multiple of an element of
$\cd^{s,\,1/16}_\ball(\wz v)$,  which together with Lemma \ref{l2.1}  implies that
$\wz v\in\dot M^{s,\,n/s}_\ball(\boz)$   and
$\|\wz v\|_{\dot M^{s,\,n/s}_\ball(\boz)}\ls
 \|h\|_{L^{n/s}(\boz)}\ls\lf( {N_0|x_1-x_2|}/L\r)^{n/s}$.
Since $\wz v\in \Delta_s(F_1,\,F_2,\,\boz)$,  we then have
$$1\ls\ccap_{\dot M^{s,\,n/s}_\ball}(F_1,\,F_2,\,\boz)\ls \lf(\frac{N_0r}L\r)^{n/s},$$
and thus $L\ls N_0|x_1-x_2|$.
So we can find a rectifiable curve $\gz^{(1)}$ joining $x_1,\,x_2$.

Now we prove the claim that $h$ is a constant multiple of an element of
$\cd^{s,\,1/16}_\ball(v)$. To this end,  we only need to chace that for $x,\,y\in\boz$
with $|x-y|\le d(x,\,\boz^\complement)/16$,
\begin{equation}\label{e3.4}
 |u(x)-u(y)|\ls|x-y|^s[h(x)+h(y)].
\end{equation}
If  $x,\,y$ lie in the same component of $\boz\setminus \overline{B(x_1,\,N_0|x_1-x_2|)}$,
then $u(x)=u(y)$ and thus \eqref{e3.4} holds.
Assume that $x,\,y$ lie in the different components of $\boz\setminus \overline{B(x_1,\,N_0|x_1-x_2|)}$.
Then the line segment joining $x$ and $y$ has a nonempty intersection with $B(x_1,\,N_0|x_1-x_2|)$
and assume it contains $w$. Moreover, $|u(x)-u(y)|\le |x-y|/L$.  Since
$$d(x,\,\boz^\complement)\le |x-w|+d(w,\,\boz^\complement)\le |x-y|
+d(w,\,\boz^\complement)\le d(x,\,\boz^\complement)/16+d(w,\,\boz^\complement)$$
and $d(w,\,\boz^\complement)\le 2 N_0|x_1-x_2|$
imply that $|x-y|\le d(x,\,\boz^\complement)\ls d(w,\,\boz^\complement)\ls N_0|x_1-x_2|$,
so  \eqref{e3.4} holds and thus gives the above claim.

Moreover, without loss of generality, we may assume that $\gz^{(1)}\cap F_i$ consists of a unique point,
  $x_i^{(1)}$, for $i=1,\,2$.
Let $  F_i^{(1)}$ be the component of $\gz^{(1)}\cap B(x_i,\,|x_i-x_i^{(1)}|/4) $ containing $x_i$,
and
$ E_i$  be the component  of $\gz^{(1)}\cap B(x_i^{(1)},\,|x_i-x_i^{(1)}|/4)$ containing $x_i^{(1)}$.
Then repeating the above procedure we can find a curve $\gz_i^{(2)}$
joining $F_i^{(1)}$  and $E_i$
such that $\ell(\gz_i^{(2)})\ls |x_i-x_i^{(1)}|$.  Denote by $x_i^{(2)}$ the unique point of $\gz_i^{(2)}\cap F_i$.
Then $|x_i-x_i^{(2)}|\le |x_i-x_i^{(1)}|/4\le |x_1-x_2|/4^2$,
$\gz^{(1)}\cup\gz_1^{(2)}\cup\gz_2^{(2)}$ contains a curve $\gz^{(2)}$ joining $x_1^{(2)}$ and $x_2^{(2)}$
with $$\ell (\gz^{(2)})\ls |x_1-x_2|+\sum_{i=1,\,2}|x_i-x_i^{(1)}|/2\ls [1+1/2]|x_1-x_2|.$$
Repeating this procedure   $k$ times until $|x_i-x_i^{(k)}|<d(x_i,\,\boz^\complement)$ for $i=1,\,2$,
we obtain a curve $\gz^{(k)}\subset\boz$
joining $x_1^{(k)}$ and $x_2^{(k)}$ with
$$\ell (\gz^{(k)})\ls |x_1-x_2|+\sum_{j=1}^k\sum_{i=1,\,2}|x_i-x_i^{(j)}|/2\ls |x_1-x_2| \sum_{j=0}^k(1+2\times4^{-j})\ls |x_1-x_2|.$$
Let $\gz$ be the union of $\gz^{(k)}$, the line segment  joining
$x_1$ and $x_1^{(k)}$  and the line segment joining $x_2$  and $x_2^{(k)}$.
Then we know $\gz \subset\boz$ joins  $x_1$ and $x_2$, and $\ell(\gz)\ls|x_1-x_2|$,
which is as desired and thus completes the proof of Theorem \ref{t1.1}.
\end{proof}

\section{Proofs of Theorems \ref{t1.2} and \ref{t1.3}\label{s4}}

The proofs of Theorems \ref{t1.2} and \ref{t1.3} consist of a sequence of auxiliary conclusions,
in particular, Theorem  \ref{t4.1}, Theorem \ref{t4.2} and Proposition \ref{p4.1} below.

We begin with several equivalent characterizations of Haj\l asz-Sobolev imbeddings,
whose proof borrows some ideas from \cite{kr93,k98,h03,hkt08b}.
In what follows, for $R\in(0,\,\fz)$ and $u\in L^1_\loc(\boz)$,
we define the maximal function $\cm^\boz_R(u)(x)$ for all $x\in\boz$  by  $$
\cm^\boz_R(u)(x)\equiv\sup_{r\in(0,\,R)}\frac1{|B(x,\,r)\cap\boz|}\int_{B(x,\,r)\cap\boz}|u(z)|\,dz.$$

\begin{thm}\label{t4.1}
Let $s\in(0,\,1]$ and $p\in(n/s,\,\fz)$.
Then the following are equivalent:

(i) $\boz$ is an $M^{s,\,p}_\ball $-imbedding domain;

(ii) $\boz$ supports the imbedding from $M^{s,\,p}_\ball(\boz)$  to $\dot M^{s-n/p,\,\fz}(\boz)$;

 (iii)  There exist positive constants $\dz_1$ and $C$  such that
for all $u\in M^{s,\,p}_\ball(\boz)$ and almost all $x,\,y\in\boz$ with $|x-y|\le \dz_1$,
$$
|u(x)-u(y)|\le C|x-y|^{s-n/p}\|u\|_{M^{s,\,p}_\ball(\boz)};
$$

 (iv)  There exist positive constants $\dz_2$ and $C$  such that
for all $u\in \dot M^{s,\,p}_\ball(\boz)$ and almost all $x,\,y\in\boz$ with $|x-y|\le \dz_2$,
\begin{equation}\label{e4.1}
|u(x)-u(y)|\le C|x-y|^{s-n/p}\|u\|_{\dot M^{s,\,p}_\ball(\boz)};
\end{equation}

 (v)  There exist positive constants $\dz_3$, $N_1$ and $C$  such that
for all $u\in \dot M^{s,\,p}_\ball(\boz)$, $g\in\cd^s_\ball(u)$ and almost all $x,\,y\in\boz$ with $|x-y|\le \dz_3$,
\begin{equation}\label{e4.2}
 |u(x)-u(y)|\le C|x-y|^{s-n/p}\lf\{\int_{B(x,\,N_1|x-y|)\cap\boz}[g(z)]^p\,dz\r\}^{1/p};
\end{equation}

 (vi)  $\boz$ is regular and there exist positive constants $\dz_4$,  $N_2$ and $C$  such that
for all $u\in \dot M^{s,\,p}_\ball(\boz)$, $g\in\cd^s_\ball(u)$ and almost all $x,\,y\in\boz$ with $|x-y|\le \dz_4$,
\begin{equation}\label{e4.3}
|u(x)-u(y)|\le C|x-y|^{s}\lf\{\cm^\boz_{N_2|x-y|}(g^p)(x)+\cm^\boz_{N_2|x-y|}(g^p)(y)\r\}^{1/p}.
\end{equation}

Moreover, if $\boz$ is bounded, then it is an $\dot M^{s,\,p}_\ball$-imbedding domain
if and only if one/all of (i) through (vi) holds.
\end{thm}

\begin{proof}[Proof of Theorem \ref{t4.1}]
We first notice that if $\boz$ is bounded, then (iv) means that $\boz$ is an $\dot M^{s,\,p}_\ball$-imbedding domain.
So it suffices to prove the equivalence of (i)  through   (vi).
Obviously,  (i)$\Rightarrow$(ii)$\Rightarrow$(iii).
In what follows, we will prove that
 (iii)$\Rightarrow$(iv)$\Rightarrow$(v)$\Rightarrow$(vi)$\Rightarrow$(iv)$\Rightarrow$(i).

(iii)$\Rightarrow $(iv). Let $u\in\dot M^{s,\,p}_\ball(\boz)$
and $g\in\cd^s_\ball(u)$ with $\|g\|_{L^p(\boz)}\ls \|u\|_{\dot M^{s,\,p}_\ball(\boz)}$.
Let  $x,\,y\in\boz$ be any pair of points satisfying $|x-y|<\dz_1$.
Without loss of generality,  we may assume that $u(y)<u(x)$
and $u(y)\le u(z)\le u(x)$ for all $z\in\boz$.
In fact, set
$$v(z)=\lf\{\begin{array}{ll}u(x) &if \ u(z)>u(x);\\
u(z) \quad& if\ u(y)\le u(z)\le u(x);\\
u(y)&if \ u(z)<u(y).
\end{array}\r. $$
Then $\cd_\ball^s(u)\subset\cd^s_\ball(v)$ and thus
$ \|v\|_{\dot M^{s,\,p}_\ball(\boz)}\ls
\|u\|_{\dot M^{s,\,p}_\ball(\boz)}$. So it suffices to prove \eqref{e4.1} for $v$.
Moreover, we may assume that $u(y)\ge 0$.
In fact, if $u(x)\le0$, then we only need to consider $-u$.
If $u(y)<0\le u(x)$, then let $u_1 =u \chi_{\{z\in\boz:\ u(z)\ge0\}}$
and $u_2=-u \chi_{\{z\in\boz:\ u(z)\le0\}}$.
Notice that
$\cd_\ball^s(u)\subset\cd^s_\ball(u_1)\cap\cd^s_\ball(u_2)$, which implies that
 $\|u_1\|_{\dot M^{s,\,p}_\ball(\boz)}+\|u_2\|_{\dot M^{s,\,p}_\ball(\boz)}\ls
\|u\|_{\dot M^{s,\,p}_\ball(\boz)}$,
and $|u(x)-u(y)|\le |u_1(x)-u_1(y)|+|u_2(x)-u_2(y)|$.
So we only need to prove \eqref{e4.1} for $u_1$ and $u_2$.

Let $\vz$ be a smooth function satisfying that
$\vz(z)=1$ for $z\in B(x,\,\dz_1)$, ${\rm supp\,} \vz\subset B(x,\,10\dz_1)$,
$0\le \vz(z)\le1$ and $|\nabla \vz(z)|\le 100/\dz_1$ for all $z\in\rn$.
Define $v(z)\equiv[u(z)-u(y)]\vz(z)$ for all $z\in\boz$.
Then it is easy to check that $v\in M^{s,\,p}_\ball(\boz)$
and $$g\vz+100(\dz_1)^{-s}[u(x)-u(y)]\chi_{\{B(x,\,10\dz_1)\cap\boz\}}\in\cd^s_\ball(v),$$
which together with $ u(y)\le u(z)\le u(x)$ for all $z\in\boz$ implies that
$$\|v\|_{ M^{s,\,p}_\ball(\boz)} \ls\|g\vz\|_{L^p(\boz)}+\|v\|_{L^{p}(\boz)}\ls
\|u\|_{M^{s,\,p}_\ball(\boz)}+[(\dz_1)^{-s+n/p}+(\dz_1)^{n/p}][u(x)-u(y)].$$
Thus, by (iii),
\begin{eqnarray*}
|u(x)-u(y)|&&=|v(x)-v(y)|\\
&&\ls|x-y|^{s-n/p}\|v\|_{ M^{s,\,p}_\ball(\boz)}\\
&&\ls|x-y|^{s-n/p}\|u\|_{\dot M^{s,\,p}_\ball(\boz)}+|x-y|^{s-n/p}[(\dz_1)^{-s+n/p}+(\dz_1)^{n/p}][u(x)-u(y)],
\end{eqnarray*}
which together with $s-n/p>0$
implies that there exists a positive constant
$\dz_2\in(0,\,\dz_1)$ such that \eqref{e4.1} holds when $|x-y|\le\dz_2$.

(iv)$\Rightarrow $(v).
Let $u\in\dot M^{s,\,p}_\ball(\boz)$ and $g\in\cd^s_\ball(u)$ with $\|g\|_{L^p(\boz)}\ls \|u\|_{\dot M^{s,\,p}_\ball(\boz)}$.
Let  $x,\,y\in\boz$ be a pair of points satisfying $|x-y|\le\dz_2$.
By an argument similar to the above,
we may assume that $0\le u(y)<u(x)$
and $u(y)\le u(z)\le u(x)$ for all $z\in\boz$.

Let $N\ge 2$ and $\vz$ be a smooth function satisfying that
$\vz(z)=1$ for $z\in B(x,\,|x-y|)$, ${\rm supp\,} \vz\subset B(x,\,N|x-y|)$,
$0\le \vz(z)\le1$ and $|\nabla \vz(z)|\le 10/(N|x-y|)$ for all $z\in\rn$.
Define $v(z)\equiv[u(z)-u(y)]\vz(z)$ for all $z\in\boz$. Then
it is easy to check that $v\in\dot M^{s,\,p}_\ball(\boz)$
and $g\vz+10(N|x-y|)^{-s}[u(x)-u(y)]\chi_{\{B(x,\,N|x-y|)\cap\boz\}}\in\cd^s_\ball(v)$,
which implies that
$\|v\|_{\dot M^{s,\,p}_\ball(\boz)} \ls
\|u\|_{\dot M^{s,\,p}_\ball(\boz)}+(N|x-y|)^{-s+n/p}[u(x)-u(y)]$.
Thus
\begin{eqnarray*}
|u(x)-u(y)|&&=|v(x)-v(y)|\\
&&\ls|x-y|^{s-n/p}\|v\|_{ \dot M^{s,\,p}_\ball(\boz)}\\
&&\ls|x-y|^{s-n/p}\lf(\int_{B(x,\,N|x-y|)\cap\boz}[g(z)]^p\,dz\r)^{1/p}+N^{-s+n/p}[u(x)-u(y)],
\end{eqnarray*}
from which and $s-n/p>0$, it follows that
there exists a positive constant $N_1\equiv N$ large enough such that if $|x-y|\le\dz_3\equiv\dz_2$, then \eqref{e4.2} holds.

(v)$\Rightarrow $(vi). We first prove that $\boz$ is regular.
For fixed $x_0\in\overline\boz$ and $0<r<\dz_3$, we define
$u(z)\equiv \frac1r d(z,\,B(x_0,\,r)^\complement)$
and $g(z)\equiv r^{-s}\chi_{\boz\cap\overline{B(x_0,\,r)}}(z)$
for all $z\in\boz$.
Then similarly to the proof of Theorem \ref{t1.1},
it is easy to check that there exists a positive constant $C$ independent of $u,\, x_0,\,r$ such that
$Cg$ is an element of $\cd^s_\ball(u)$, which implies that
$u\in\dot M^{s,\,p}_\ball(\boz)$
and $\|u\|_{\dot M^{s,\,p}_\ball(\boz)}\ls r^{-s}|\boz\cap B(x_0,\,r)|^{1/p}$. By this and \eqref{e4.2}, we further have
$$1\le r^{s-n/p} r^{-s}|B(x_0,\,r)\cap\boz|^{1/p}$$
and thus $|B(x_0,\,r)\cap\boz|\ge r^n$.

Now, let $u\in\dot M^{s,\,p}_\ball(\boz)$ and $g\in\cd^s_\ball(u)$ with $\|g\|_{L^p(\boz)}\ls \|u\|_{\dot M^{s,\,p}_\ball(\boz)}$.
Let  $x,\,y\in\boz$ be a pair of points satisfying $|x-y|<\dz_3/10$.
Since $\boz$ is regular, we then have
\begin{eqnarray*}
 |u(x)-u_{B(x,\,2|x-y|)\cap\boz}|&&\le\sum_{j=-1}^\fz|u_{B(x,\,2^{-j}|x-y|)\cap\boz}-u_{B(x,\,2^{-j-1}|x-y|)\cap\boz}|\\
&&\ls\sum_{j=-1}^\fz (2^{-j}|x-y|)^{s-n/p}\lf(\int_{B(x,\,2^{-j}N_1|x-y|)\cap\boz}[g(z)]^p\,dz\r)^{1/p}\\
&&\ls|x-y|^{s-n/p}\lf\{\cm^\boz_{2N_1|x-y|} (g^p)\r\}^{1/p}(x).
\end{eqnarray*}
Similarly, we can prove that
\begin{eqnarray*}
 |u(y)-u_{B(x,\,2|x-y|)\cap\boz}|
&&\ls|x-y|^{s-n/p}\lf\{\cm^\boz_{2N_1|x-y|} (g^p)\r\}^{1/p}(y).
\end{eqnarray*}
Thus we obtain \eqref{e4.3} with $N_2\equiv 2N_1$ and $\dz_4\equiv \dz_3/10$ and hence (vi).

(vi)$\Rightarrow $(iv).
Let $u\in\dot M^{s,\,p}_\ball(\boz)$ and $g\in\cd^s_\ball(u)$
with $\|g\|_{L^p(\boz)}\ls \|u\|_{\dot M^{s,\,p}_\ball(\boz)}$.
By a slight modification of the
proof of \cite[Theorem 9.5]{h03},
we know that \eqref{e4.3} implies that for all $r\le \dz_4/2$,
$$\bint_{B(x,\,r)\cap\boz}|u(z)-u_{B(x,\,r)\cap\boz}|\,dz\ls r^s\lf(\bint_{B(x,\,6Nr)\cap\boz}[g(z)]^p\,dz\r)^{1/p},$$
where $N$ is a positive constant independent of $u$, $g$ and $r$.
By this and the fact that $\boz$ is regular, for $|x-y|\le \dz_4/4$, we have
\begin{eqnarray*}
 |u(x)-u_{B(x,\,2|x-y|)\cap\boz}|&&\le\sum_{j=-1}^\fz|u_{B(x,\,2^{-j}|x-y|)\cap\boz}-u_{B(x,\,2^{-j-1}|x-y|)\cap\boz}|\\
&&\ls\sum_{j=-1}^\fz (2^{-j}|x-y|)^{s-n/p}\lf(\int_{B(x,\,6\cdot2^{-j}N|x-y|)\cap\boz}[g(z)]^p\,dz\r)^{1/p}\\
&&\ls|x-y|^{s-n/p}\lf(\int_{\boz}[g(z)]^p\,dz\r)^{1/p},
\end{eqnarray*}
which also holds for $|u(y)-u_{B(x,\,2|x-y|)\cap\boz}|$ by a similar argument.
Thus  we have \eqref{e4.1} and (iv).

(iv)$\Rightarrow $(i).
By an argument similar to that used in (v)$\Rightarrow $(vi), we know that $\boz$ is regular.
Let $u\in M^{s,\,p}_\ball(\boz)$. Then we only need to prove that $\|u\|_{L^\fz(\boz)}\ls\|u\|_{M^{s,\,p}_\ball(\boz)}$
and for almost all $x,\,y\in\boz$, $|u(x)-u(y)|\ls \|u\|_{M^{s,\,p}_\ball(\boz)}$.
In fact, for almost $x\in\boz$, by the H\"older inequality and \eqref{e4.1}, we have
$$|u(x)|\le |u(x)-u_{B(x,\,\dz_2/2)\cap\boz}|+\|u\|_{L^p(\boz)}\ls \|u\|_{M^{s,\,p}_\ball(\boz)},$$
which implies that $\|u\|_{L^\fz(\boz)}\ls\|u\|_{M^{s,\,p}_\ball(\boz)}$.
Moreover, for almost all $ x,\,y\in\boz$,  if $|x-y|\ge\dz_2/2$, then
$|u(x)-u(y)|\ls \|u\|_{L^\fz(\boz)}\ls\|u\|_{M^{s,\,p}_\ball(\boz)}$;
if $|x-y|<\dz_2/2$, then \eqref{e4.1} yields that
$|u(x)-u(y)|\ls  \|u\|_{M^{s,\,p}_\ball(\boz)}$. This shows (i) and hence  finishes the proof of Theorem \ref{t4.2}.
\end{proof}

\begin{rem}\label{r4.1}\rm
We point out that if $\boz$ is a bounded $\dot M^{s,\,p}_\ball$-imbedding domain,
then  (iv) holds with $\dz_2=\diam\boz$ and hence, by the proofs of (iv)$\Rightarrow $(v) and (iv)$\Rightarrow $(v),
we can further take $\dz_3= \diam\boz$ in (v) and also $\dz_4=\diam\boz$ in (vi).
\end{rem}

By an argument similar to the proofs of (iv)$\Rightarrow$(v)$\Rightarrow$(vi)
in Theorem \ref{t4.1} and the observation as in Remark \ref{r4.1}, we have the following conclusion.
\begin{cor}\label{c4.1}
Let $s\in(0,\,1]$ and $p,\,\wz p\in(n/s,\,\fz)$ with $p<\wz p$.
Assume that there exists a positive constant $C$ such that for all
$u\in\dot M^{s,\,\wz p}_\ball(\boz)$, $\eqref{e4.1}$ of Theorem \ref{t4.1} holds with the same constants, namely,
for almost all $x,\,y\in\boz$ with $|x-y|\le \dz_2$,
$$|u(x)-u(y)|\le C|x-y|^{s-n/p}\|u\|_{\dot M^{s,\,p}_\ball(\boz)}.$$
Then for all $u\in\dot M^{s,\,\wz p}_\ball(\boz)$ and $g\in \cd^s_\ball(u)$,
 \eqref{e4.2} and \eqref{e4.3} of Theorem \ref{t4.1} still hold with the same constants.
Moreover, if $\boz$ is a bounded domain, then \eqref{e4.2} and \eqref{e4.3} holds with $\dz_3=\dz_4=\diam\boz$.
\end{cor}

The following conclusion clarifies the relations between Haj\l asz-Sobolev extensions
and imbeddings to some extent, and hence
generalizes  \cite[Theorem A]{k98}.

\begin{thm}\label{t4.2}
If $\boz$ is an  $M^{s,\,p}_\ball $-imbedding (resp. a bounded  $\dot M^{s,\,p}_\ball $-imbedding) domain for
some $s\in(0,\,1]$ and $p\in(n/s,\,\fz)$, then it is an $M^{t,\,q}_\ball$-extension (resp. $\dot M^{t,\,q}_\ball$-extension) domain
for all $t\in[s,\,1]$ and $q\in(n/t,\,\fz)$ satisfying
$t-n/q>s-n/p$.
\end{thm}

To prove Theorem \ref{t4.2}, we need the following conclusion, which is essentially established
in \cite{hkt08b} and also \cite{hkt08,s06}.

\begin{lem}\label{l4.1}
Let $s\in(0,\,1]$ and $p\in(1,\,\fz)$.
A (bounded) domain $\boz\subset\rn$ is an  $M^{s,\,p}$-extension (resp. $\dot M^{s,\,p}$-extension) domain if and only if
$\boz$ is regular.
A  (bounded) domain $\boz\subset\rn$ is an $M_\ball^{s,\,p}$-extension (resp. $\dot M^{s,\,p}_\ball$-extension)
domain if and only if
$\boz$ is regular and $\dot M_\ball^{s,\,p}(\boz)= \dot M^{s,\,p}(\boz)$.  Moreover, the extension operators can be assumed to be linear.
\end{lem}

Notice that for every $s\in(0,\,1]$,
$(\rn,\,d_s,\,dx)$  is an
Ahlfors $n/s$-regular metric measure space,
where $d_s(x,\,y)=|x-y|^s$ and $dx$ denotes the Lebesgue measure.
Since $  M^{s,\,p}(\boz)=  M^{1,\,p}(\boz,\,d_s,\,dx)$,
Lemma \ref{l4.1} for the inhomogeneous Haj\l asz-Sobolev spaces is given by
\cite[Theorem 5]{hkt08b} (see also \cite[Theorem 1.3]{s06}).
When $\boz$ is bounded, by an argument similar to that of \cite[Theorem 5]{hkt08b} (and also \cite[Theorem 1.3]{s06}),
Lemma \ref{l4.1} still holds for the homogeneous Haj\l asz-Sobolev spaces.

\begin{proof}[Proof of Theorem \ref{t4.2}]
First, we point out that it suffices to prove Theorem \ref{t4.2} for inhomogeneous Haj\l asz-Sobolev spaces.
To see this, assume that $\boz$ is a bounded $\dot M^{s,\,p}_\ball$-imbedding domain.
By Theorem \ref{t4.1}, then $\boz$ is an $M^{s,\,p}_\ball$-imbedding domain.
If Theorem \ref{t4.2} holds for the inhomogeneous Haj\l asz-Sobolev spaces, then
for all $t\in[s,\,1]$ and $q\in(n/t,\,\fz)$ satisfying $t-n/q>p-n/s$,
$\boz$ is an $M^{t,\,q}_\ball$-extension domain and hence
an $M^{t,\,q}_\ball$-imbedding domain,
which together with  Theorem \ref{t4.1} again implies that
$\boz$ is an $\dot M^{t,\,q}_\ball$-imbedding domain.
So, it further suffices to show that
 a bounded $\dot M^{s,\,p}_\ball$-imbedding domain
is an $\dot M^{s,\,q}_\ball$-extension domain for all $q\in(p,\,\fz)$.
To this end, let $u\in M^{s,\,q}_\ball(\boz)$ and $g\in\cd^s_\ball(u)$ with $\|g\|_{L^q(\boz)}\ls\|u\|_{\dot M^{s,\,q}_\ball(\boz)}$.
By Corollary \ref{c4.1}, we know that \eqref{e4.3} holds for all $x,\,y\in\boz$,
which  means that $\lf\{\cm^\boz_{N_2\dz_4}(g^p)\r\}^{1/p}$ is a constant multiple of an element of $\cd^s(u)$.
Hence by the $L^{q/p}(\boz)$-boundedness of $\cm^\boz_{N_2\dz_4}$, we have
$u\in M^{s,\,q} (\boz)$ and
$$\|u\|_{\dot M^{s,\,q} (\boz)}\ls \|\lf\{\cm^\boz_{N_2\dz_4}(g^p)\r\}^{1/p}\|_{L^q(\boz)}\ls
\|g\|_{L^q(\boz)}
\ls\|u\|_{\dot M^{s,\,q}_\ball(\boz)}.$$
By Lemma \ref{l4.1}, we deduce that
$\boz$ is an $\dot M^{s,\,q}_\ball$-extension domain for all $q>p$.

To prove Theorem \ref{t4.2} for the inhomogeneous Haj\l asz-Sobolev spaces,
assume that $\boz$ is an $M^{s,\,p}_\ball$-imbedding domain.
Notice that, by Theorem \ref{t4.1},
$\boz$ is regular.
So by Lemma \ref{l4.1}, it suffices to prove
$M^{t,\,q}_\ball(\boz)=M^{t,\,q}(\boz)$.
Obviously,  $M^{t,\,q}(\boz)\subset M^{t,\,q}_\ball(\boz)$,
so we only need to
prove that $M^{t,\,q}_\ball(\boz)\subset M^{t,\,q}(\boz)$.
We consider the following two cases.

  {\it Case $t=s$.}
Let $u\in M^{s,\,q}_\ball(\boz)$ and $g\in\cd^s_\ball(u)$ with $\|g\|_{L^q(\boz)}\ls\|u\|_{\dot M^{s,\,q}_\ball(\boz)}$.
Let $\dz\equiv\min\{\dz_3,\,\dz_4\}$ and $N\equiv\max\{N_1,\,N_2\}$,
where $N_1,\,N_2,\,\dz_1,\,\dz_2$ are as in Theorem \ref{t4.1}.
 Then for all $x,\,y\in\boz$ with $|x-y|\le \dz$,  by Corollary \ref{c4.1},
we know that   \eqref{e4.3}  holds.
Notice that by Theorem \ref{t4.1}, $\boz$ is regular,
and by Corollary \ref{c4.1} again, \eqref{e4.2} also holds. So for almost all $x\in\boz$,
we have
\begin{eqnarray*}
|u(x)|&&\le |u(x)-u_{B(x,\,\dz )\cap\boz}|+ |u_{B(x,\,\dz )\cap\boz}|\\
&&\ls \bint_{B(x,\,N \dz )\cap\boz}|g(z)|\,dz+\bint_{B(x,\,\dz )\cap\boz}|u(z)|\,dz
 \\
&&\ls [\cm^\boz_{N\dz }(g^p)(x)]^{1/p}+ \cm^\boz_{N\dz }(u)(x).
\end{eqnarray*}
which implies that, for all $x,\,y\in\boz$ with  $|x-y|\ge \dz$,
$$|u(x)-u(y)|\ls |x-y|^{s}\{\cm^\boz_{N\dz }(u)(x)+[\cm^\boz_{N\dz }(g^p)(x)]^{1/p}+
\cm^\boz_{N\dz }(u)(y)+[\cm^\boz_{N\dz }(g^p)(y)]^{1/p}\}.$$
Therefore,
$\cm^\boz_{N \dz }(u) +\{\cm^\boz_{ N  \dz }(g^p)\}^{1/p}$ is a constant multiple of an element of $\cd^s(u)$,
which together with $q>p$ and
the $L^{q/p}(\rn)$-boundedness of $\cm$  implies that $u\in M^{s,\,q}(\boz)$ and $\|u\|_{M^{s,\,q}(\boz)}\ls\|u\|_{M_\ball^{s,\,q}(\boz)}$.

{\it Case $t>s$.} By the  {\it Case $t=s$} and the conclusion of its proof, it suffices to prove that $\boz$ is an
$M^{t,\,q}_\ball$-imbedding domain for all $q\in(n/t,\,\fz)$ satisfying $q-n/t>p-n/s$.
Let $u\in M^{t,\,q}_\ball(\boz)$ and $g\in\cd^t_\ball(u)$ with $\|g\|_{L^q(\boz)}\ls\|u\|_{\dot M^{t,\,q}_\ball(\boz)}$.
Let $x,\,y\in\boz$ with $|x-y|\le\min\{\dz_2,\,\dz_3\}$, where $\dz_2$ and $\dz_3$ are as in Theorem \ref{t4.1}.
Then by Theorem \ref{t4.1}, it further suffices to prove that there exists a positive constant $C$ independent of $u$,
$x$ and $y$ such that
\begin{equation}\label{e4.4}
|u(x)-u(y)|\ls
|x-y|^{t-n/q}\lf(\int_{B(x,\,N|x-y|)\cap\boz}[g(z)]^{q}\,dz\r)^{1/q}.
\end{equation}

Without loss of generality, we assume that $0\le u(y)<u(x)$.
Let $\vz$ be a smooth function satisfying that
$\vz(z)=1$ for $z\in B(x,\,|x-y|)$, ${\rm supp\,} \vz\subset B(x,\,N|x-y|)$,
$0\le \vz(z)\le1$ and $|\nabla \vz(z)|\le 10/(N|x-y|)$ for all $z\in\rn$.
Define $v(z)\equiv[u(z)-u(y)]\vz(z)$ for all $z\in\boz$.
Then it is easy to see that $$h\equiv g\vz+10(N|x-y|)^{-t}[u(x)-u(y)]\chi_{\{B(x,\,N|x-y|)\cap\boz\}}\in\cd^t_\ball(v),$$
which implies that $v\in\dot M^{t,\,q}_\ball(\boz)$
and
\begin{equation}\label{e4.5}
 \|v\|_{\dot M^{t,\,q}_\ball(\boz)} \ls\|h\|_{L^q(\boz)}\ls
\|g\vz\|_{L^q(\boz)}+(N|x-y|)^{-t+n/q}[u(x)-u(y)].
\end{equation}
Moreover, we claim that
$\wz h\equiv \cm^{(t-s)}(h)$ is a constant multiple of an element of $\cd^{s,\,1/8}_\ball(v)$,
where
 $$\cm^{(t-s)}(h)(z)\equiv\sup_{0<\wz r<d(z,\,\boz^\complement)/2} (\wz r)^{t-s} \bint_{B(z,\,\wz r)}|h(w)|\,dw.$$
In fact, for every pair of points  $z,\,w\in\boz$ and $0<|z-w|< \frac18d(z,\,\partial\boz)$,
\begin{eqnarray*}
|v(z)-v_{B(z,\,2|z-w|)}|&&= \sum_{j=-1}^\fz|v_{B(z,\,2^{-j}|z-w|)}-v_{B(z,\,2^{-j-1}|z-w|)}|\\
&&\ls\sum_{j=-1}^\fz(2^{-j}|z-w|)^{t-n}\int_{B(z,\,2^{-j }|z-w|)}|h(\wz z)|\,d\wz z\\
&&\ls |z-w|^{s}\cm^{(t-s)}(h)(z),
\end{eqnarray*}
and similarly,  $|v(w)-v_{B(z,\,2|z-w|)}|\ls |z-w|^{s}\cm^{(t-s)}(h)(z)$, which imply the above claim.
Let $\wz q\in(p,\,\fz)$ such that  $1/q-1/\wz q=(t-s)/n$.
 Then by the boundedness from $L^{ q}(\rn)$ to $L^{\wz q}(\rn)$ of $\cm^{(t-s)}$ (see \cite{s93}),
we have
$
\|\cm^{(t-s)}(h)\|_{L^{\wz q}(\boz)}\ls \|h\|_{L^q(\boz)},
$
which together with Lemma \ref{l2.1} implies that
$v\in  \dot M^{s,\,\wz q}_\ball(\boz)$  and
$\|v\|_{\dot M^{s,\,\wz q}_\ball(\boz)}\ls   \|h\|_{L^q(\boz)}$.
Since $\wz q>p$, then by  {\it Case $t=s$}, $\boz$ is an $M^{s,\,\wz q}_\ball$-imbedding domain.
Then by Theorem \ref{t4.1}(iv), $\|v\|_{\dot M^{s,\,\wz q}_\ball(\boz)}\ls   \|h\|_{L^q(\boz)}$,
 \eqref{e4.5} and $1/q-1/\wz q=(t-s)/n$,
we have
\begin{eqnarray*}
|u(x)-u(y)|&&=|v(x)-v(y)|  \ls|x-y|^{s-n/\wz q} \|v\|_{\dot M^{s,\,\wz q}_\ball(\boz)}\\
&&\ls |x-y|^{t-n/q}\lf(\int_{B(x,\,N|x-y|)\cap\boz}[g(z)]^{q}\,dz\r)^{1/q}+N^{-t+n/q}[u(x)-u(y)].
\end{eqnarray*}
If $N$ is large enough, by $t-n/q>0$, we then have \eqref{e4.4}, which completes the proof of Theorem \ref{t4.2}.
\end{proof}

\begin{rem}\label{r4.2}\rm
(i) Recall that if $\boz$ is an  $M^{s,\,p}_\ball$-extension
(resp. a bounded $\dot M^{s,\,p}_\ball$-extension) domain for some $s\in(0,\,1]$ and $p\in(n/s,\,\fz)$,
then it is an $M^{s,\,p}_\ball $-imbedding
(resp. $\dot M^{s,\,p}_\ball$-imbedding) domain,
and hence by Theorem \ref{t4.2},
$\boz$ is an $M^{t,\,q}_\ball $-extension/-imbedding
(resp. $\dot M^{t,\,q}_\ball$-extension/-imbedding) domain for all $t\in[s,\,1]$ and $q\in(n/t,\,\fz)$ satisfying
$t-n/q> s-n/p$.

(ii) We also point out that the extension operators in Theorem \ref{t4.2} can be assumed to be linear due to Lemma \ref{l4.1} below.
\end{rem}

By checking the {\it Case $t=s$} in the proof of Theorem \ref{t4.2},
we obtain a weak version of Theorem \ref{t4.2} when $t=s$
with the aid of Corollary \ref{c4.1}, which
reduces the assumption that $\boz$ is an $\dot M^{s,\,p}$-imbedding appeared in Theorem \ref{t4.2} slightly
and plays an important role in the proof of  Theorem \ref{t1.2}.

\begin{cor}\label{c4.2}
Let $s\in(0,\,1]$, $p,\,\wz p\in(n/s,\,\fz)$ with $p<\wz p$ and $\boz$ be a bounded domain.
Suppose that there exists a positive constant $C$ such that if
$u\in\dot M^{s,\,\wz p}_\ball(\boz)$, then $u\in \dot M^{s-  n/p,\,\fz}(\boz)$ and $
 \|u\|_{\dot M^{s- n/p,\,\fz}(\boz)}\le C\|u\|_{\dot M^{s,\,  p}_\ball(\boz)}.$
Then $\boz$ is an $\dot M^{s,\,\wz p}_\ball$-extension domain.
\end{cor}

The following weak self-improving property established
by Shvartsman \cite[Theorem 1.4]{s10} plays an crucial role in the proof of Theorem \ref{t1.2}.

\begin{prop}\label{p4.1}
Let $\az\in(0,\,1)$ and   $\boz$ be a weak $\az$-cigar domain in $\rn$.
There exist a  constant  $\az^\ast\in(0,\,\az)$, and constants $\tz,\,C>0$ such that the following are true:

For every $\ez>0$ and $x,\,y\in\boz$ with $|x-y|\le\tz$, there exist a rectifiable curve $\Gamma\subset\boz$
joining $x$ and $y$ and a subset $\wz \Gamma\subset \Gamma$ consisting of a finite number of arcs such that the following
conditions are satisfied:

(i)  For every $\tau\in[\az^\ast,\,\az]$,
\begin{equation}\label{e4.6}
 \int_{\wz\Gamma}[d(z,\,\boz^\complement)]^{\tau-1}\,|dz|\le C|x-y|^\tau.
\end{equation}
In addition, for every ball $B$ centered in $\wz \Gamma$ of radius at most $|x-y|$,
$\diam B\le C\ell(B\cap \wz \Gamma)$.

(ii) We have $\ell(\Gamma)\le C|x-y|$ and
\begin{equation}\label{e4.7}
\ell(\Gamma\setminus\wz\Gamma)<\ez.
\end{equation}
Moreover,
\begin{equation}\label{e4.8}
 \int_{\Gamma\setminus\wz\Gamma}[d(z,\,\boz^\complement)]^{\az-1}\,|dz|\le C|x-y|^\az.
\end{equation}
\end{prop}

\begin{proof}[Proof of Theorem \ref{t1.2}]
(i) Assume that $\boz$ is a weak $\az$-cigar domain with $\az\in(0,\,1)$.
Let $s\in(0,\,1]$  and   $\wz p\equiv (n-\az)/(s-\az) $ be fixed.
Let $\az^\ast\in(0,\,\az)$, $\tz$ and $C_1$ be as in Proposition \ref{p4.1},
and also $p^\ast\equiv(n-\az^\ast)/(s-\az^\ast) $.
Then   $\az^\ast=(p^\ast s-n)/(p^\ast-1)$,
$\az =(\wz p  s-n)/(\wz p -1)$ and $n/s<p^\ast<  \wz p<\fz $.
Moreover, without loss of generality, we may assume that $\tz=\diam\boz$.

We first claim that for every $ p\in[p^\ast,\, \wz p]$, there exists a positive constant $C$ such that for all
$u\in M^{s,\, \wz p}_\ball(\boz)$,
$u\in \dot M^{s-  n/p,\,\fz}(\boz)$  and
\begin{equation}\label{e4.9}
 \|u\|_{\dot M^{s- n/p,\,\fz}(\boz)}\le C\|u\|_{\dot M^{s,\,  p}_\ball(\boz)}.
\end{equation}

Then Theorem \ref{t1.2}(i) then follows from this claim.
Indeed, the above claim implies that
 $\boz$ is an $\dot M^{s,\,\wz p}_\ball$-imbedding domain,
which together with Theorem \ref{t4.2} further yields that
$\boz$ is an $\dot M^{s,\,p}_\ball$-extension and $\dot M^{s,\,p}_\ball$-imbedding domain
 for all $p\in(\wz p,\,\fz)$.
Moreover, the above claim together with Corollary \ref{c4.2} implies that $\boz$ is an
$\dot M^{s,\, \wz p}_\ball$-extension domain.

Now we turn to the proof of the above claim.
Without loss of generality, we may also assume that
$\|u\|_{\dot M^{s,\, p}_\ball(\boz)}=1$.
Then
the proof of \eqref{e4.9} is reduced to checking that for all Lebesgue points $x,\,y\in\boz$ of $u$,
\begin{equation}\label{e4.10}
 |u(x)-u(y)|\ls|x-y|^{s-n/p}.
\end{equation}

To this end, take $g\in \cd^s(u)\cap L^{p }(\boz)\cap L^{\wz p}(\boz)$ such that $1\le\|g\|_{L^{p}(\boz)}\le2$.
If $|x-y|\le d(x,\,\boz^\complement) /4$, then
\begin{eqnarray}\label{e4.11}
|u(x)-u_{B(x,\,2|x-y|)}|&&\le
\sum_{j=0}^\fz |u_{B(x,\,2^{-j}|x-y|)}-u_{B(x,\,2^{-j+1}|x-y| )}| \\
&&\ls\sum_{j=0}^\fz \bint_{B(x,\,2^{-j }|x-y|)}|u(z)-u_{B(x,\,2^{-j}|x-y| )}|\,dz\nonumber\\
&&\ls \sum_{j=-1}^\fz [2^{-j}|x-y|]^s\bint_{B(x,\,2^{-j}|x-y|)} g(z)\,dz\nonumber\\
&&\ls \sum_{j=-1}^\fz  [2^{-j}|x-y|]^{s-n/p}
\lf\{\int_{B(x,\,2^{-j }|x-y|)} [g(z)]^{p}\,dz\r\}^{1/p}\nonumber\\
&&\ls    |x-y|^{s-n/p}
\lf\{\int_{B(x,\, 2|x-y|)} [g(z)]^{p}\,dz\r\}^{1/p},\nonumber
\end{eqnarray}
which together with  $\|g\|_{L^{p}(\boz)}\sim1$ gives that
$|u(x)-u_{B(x,\,2|x-y|)}|\ls|x-y|^{s-n/p}.$  Similarly, we can prove that
$|u(y)-u_{B(x,\,2|x-y|)}|\ls|x-y|^{s-n/p}$. Thus,   $$
|u(x)-u(y)| \le |u(x)-u_{B(x,\,2|x-y|)}|+|u(y)-u_{B(x,\,2|x-y|)}|\ls|x-y|^{s-n/p},$$
which gives \eqref{e4.10} when $|x-y|\le d(x,\,\boz^\complement) /4$.

Now we assume that
$|x-y|>\max\{d(x,\,\boz^\complement), d(y,\,\boz^\complement)\}/4$.
 Then for  $\ez>0$ fixed,  let $\Gamma$ joining $x$ and $y$,
and a subset  $\wz\Gamma$ of $\Gamma$ be  as in Proposition \ref{p4.1}.
Then by \eqref{e4.6} and \eqref{e4.8}, we have
$$\int_\Gamma [d(z,\,\boz^\complement)]^{(\wz ps-n)/(\wz p-1)-1}\,|dz|\ls |x-y|^ {(\wz ps-n)/(\wz p-1)}.$$
By using the Bescovitch covering lemma (see \cite{s93}) and some arguments similar to
these in the proofs of \cite[Theorem 4.1] {bk96} and \cite[Lemma 3.2]{s10},
we can find a family of balls $\cb\equiv\{B_i\}_{i=1}^N$ such that

 a) $B_i\equiv B(z_i,\,d(z_i,\,\boz^\complement)/50)$
with $z_i\in\Gamma$ for all  $i=0,\,\cdots,\,N$, $z_0=x$ and $z_N=y$;

 b) $B_i\cap B_{i+1}\ne\emptyset$ for all $i=0,\,\cdots,\,N-1$;

 c) $\sum_{i=1}^N\chi_{2B_i}(z)\le C_2$
for all $z\in\boz$, where the constant $C_2$ only depends on the dimension $n$.

Let $w_i\in B_i\cap B_{i+1}$ for all $i=0,\,\cdots,\,N-1$.
Then by the properties a) and b) above, \eqref{e4.11}  and an argument similar to \eqref{e4.11}, we have
\begin{eqnarray*}
|u(x)-u(y)|&&\le \sum_{i=0}^{N-1}(|u(z_i)-u(w_i)|+|u(w_i)-u(z_{i+1})|)\\
&&\le \sum_{i=0}^{N-1}(|u(z_i)-u_{B(z_i,\,2|w_i-z_i|)}|+|u(w_i)-u_{B(z_i,\,2|w_i-z_i|)}|\\
&&\quad \quad+
|u(w_i)- u_{B(z_{i+1},\,2|w_i-z_ {i+1}|)}|+|u(z_{i+1})- u_{B(z_{i+1},\,2|w_i-z_ {i+1}|)}|)\\
&&\ls\sum_{i=0}^{N}[d(z_i,\,\boz^\complement)]^{s-n/p}
\lf(\int_{2B_i} [g(z)]^{p}\,dz\r)^{1/p}.
\end{eqnarray*}
Let $\wz \cb $ be the collection of all $B_i\in\cb$ such that $B_i\cap\wz\Gamma \ne\emptyset$.
Then by the H\"older inequality, we have
\begin{eqnarray*}
|u(x)-u(y)|
&&\ls\sum_{B_i\in\wz \cb  }[d(z_i,\,\boz^\complement)]^{s-n/p}
\lf(\int_{2B_i} [g(z)]^{p}\,dz\r)^{1/p}\\
&&\quad
+ \sum_{B_i\in\cb\setminus \wz\cb}[d(z_i,\,\boz^\complement)]^{s-n/\wz p}
\lf(\int_{2B_i} [g(z)]^{\wz p}\,dz\r)^{1/\wz p}\\
&&\ls  \lf(\sum_{B_i\in\wz \cb  }[d(z_i,\,\boz^\complement)]^{( p s-n)/(p-1)}\r)^{(p-1)/p}
\lf(\sum_{B_i\in \cb }\int_{2B_i} [g(z)]^{p}\,dz\r)^{1/p}\\
&&\quad+
\lf(\sum_{i=0  }^N[d(z_i,\,\boz^\complement)]^{( \wz ps -n)/(\wz p -1)}\r)^{(\wz p -1)/\wz p}
\lf(\sum_{B_i\in\cb\setminus\wz \cb  }\int_{2B_i} [g(z)]^{\wz p}\,dz\r)^{1/\wz p}\\
&&\equiv I_1+I_2.
\end{eqnarray*}

To estimate $I_1$, for each $B_i\in\wz\cb$, take $\wz z_i\in B_i\cap \wz\Gamma$.
It is easy to see that
$$\frac{24}{25} d(z_i,\,\boz^\complement)\le d(z,\,\boz^\complement)\le \frac{25}{24} d(z_i,\,\boz^\complement),$$
for all $z\in 2B_i$, especially, for $z=\wz z_i$.
By this and Proposition \ref{p4.1},
$d(\wz z_i,\,\boz^\complement)\ls \ell(\wz\Gamma\cap 2B_i)$. Thus,
$$[d(z_i,\,\boz^\complement)]^{( p  s-n)/(p -1) }\ls
\int_{\wz \Gamma\cap 2B_i} [d(z,\,\boz^\complement)]^{( p  s-n)/(p -1)-1}|dz|,$$
which together with $\sum_{B_i}\chi_{2B_i}\ls1$ and \eqref{e4.6} of Proposition \ref{p4.1} with $\tau\equiv( p s-n)/(p -1)\in[\az^\ast,\,\az]$
implies that
\begin{eqnarray*}
 I_1&&\ls\lf(\sum_{B_i\in\wz \cb  }
\int_{\wz \Gamma\cap 2B_i} [d(z,\,\boz^\complement)]^{( p s-n)/(p -1)-1}|dz|\r)^{(p-1)/p}
\lf(\sum_{B_i\in \cb }\int_{2B_i} [g(z)]^{p}\,dz\r)^{1/p}\\
&&\ls \lf( \int_{\wz \Gamma } [d(z,\,\boz^\complement)]^{( p  s-n)/(p -1)-1}|dz|\r)^{(p -1)/p }
\lf( \int_\boz [g(z)]^{p}\,dz\r)^{1/p}\\
&&\ls|x-y|^{s-n/p}.
\end{eqnarray*}
By an argument as in  the estimate of $I_1$, we then have

\begin{eqnarray*}
 I_2&&\ls\lf(\sum_{B_i\in \cb  }\int_{  \Gamma\cap 2B_i} [d(z,\,\boz^\complement)]^{(\wz ps-n)/(\wz p-1)-1}|dz|\r)^{(\wz p -1)/\wz p }
\lf(\sum_{B_i\in \cb\setminus \wz \cb }\int_{2B_i} [g(z)]^{\wz p }\,dz\r)^{1/\wz p }\\
&&\ls|x-y|^{s-n/\wz p}\lf( \int_{\cup_{B_i\in \cb\setminus \wz \cb}2B_i }  [g(z)]^{\wz p}\,dz\r)^{1/\wz p }.
\end{eqnarray*}
Moreover,  notice that if $B_i\in\cb\setminus\wz\cb$,
then $B_i\cap\Gamma\subset \Gamma\setminus  \wz\Gamma$, which together with \eqref{e4.7} implies that
$$\sum_{B_i\in\cb\setminus\wz\cb}|B_i|^{1/n}
\le  \sum_{B_i\in\cb\setminus\wz\cb}\int_{\Gamma\setminus\wz\Gamma}\chi_{B_i}(z)\,|dz|
\le \int_{\Gamma\setminus\wz\Gamma} \sum_{B_i\in\cb }\chi_{B_i}(z)\,|dz|
\ls  \ell(\Gamma\setminus\wz\Gamma)\ls\ez$$
and hence
$$|\cup_{B_i\in\cb\setminus\wz\cb}2B_i|\ls \sum_{B_i\in\cb\setminus\wz\cb}|B_i|
\ls \lf(\sum_{B_i\in\cb\setminus\wz\cb}|B_i|^{1/n} \r)^{1/n}\ls\ez.$$
Since the constants that appeared in the estimates of $I_1$ and $I_2$ are independent of $\ez$,
 by $g\in L^p(\boz)$ and absolute continuity of integral with respect to measure,
we can choose $\ez$ small enough such that $I_2\le |x-y|^{s-n/p}$.
Combining the estimates of $I_1$ and $I_2$,
we obtain \eqref{e4.10} for almost all $x,\,y\in\boz$. This finishes the proof of Theorem \ref{t1.2}(i).

(ii)  Assume that $\boz$ has the slice property as in Definition \ref{d2.4} with constant $C_0$
and also that $\boz$ is a bounded $\dot M^{s,\,p}_\ball$-imbedding domain.
Observe that, by \cite{bk96}, a bounded weak $\az$-cigar domain with $\az\in(0,\,1)$ is weak $\bz$-cigar domain
for all $\bz\in(\az,\,1]$. So it suffices to prove that $\boz$ is a weak $(ps-n)/(p-1)$-cigar domain.
Without loss of generality, we may assume that $j\ge2$ in Definition \ref{d2.4}.
Indeed, if $j=0$, then necessarily $x=y$. If $j=1$, by Definition \ref{d2.4}(iii) and (iv), we have that
 $d(z,\,\boz^\complement)\gs d(x,\,\boz^\complement)$ for all $z\in S_0\cap \gz$
and $d(z,\,\boz^\complement)\gs d(y,\,\boz^\complement)$ for all $z\in S_1\cap \gz$.
This together with Definition \ref{d2.4}(iii) implies that
\begin{eqnarray*}
 \int_\gz [d(z,\,\boz^\complement)]^{(ps-n)/(p-1)-1}\,|dz|&&\ls
[d(x,\,\boz^\complement)]^{(ps-n)/(p-1)}+[d(y,\,\boz^\complement)]^{(ps-n)/(p-1)}\\
&&\ls|x-y|^{(ps-n)/(p-1)},
\end{eqnarray*}
as desired.

Suppose  thus $x,\,y\in\boz$ are fixed, and
  $\gz$ and
$\{S_i\}_{i=0}^j$ be as in Definition \ref{d2.4} with $j\ge2$.
For each $i=1,\,\cdots,\,j-1$, define   function $u_i$ by setting
$u_i(z)\equiv\inf_{\wz\gz}\ell(\wz\gz\cap S_i)$ for all $z\in\boz$,
where the infimum is taken over all the rectifiable curves $\wz \gz$ joining $x$ and $z$.
Obviously, $u_i(z)=0$ for $z\in\cup_{k=0}^{i-1}S_k$ and $u_i(z)$ is a constant for $z\in\cup_{k=i+1}^{j}S_k$.
Similarly to the proof of Theorem \ref{t1.1},
there exists a positive constant $\wz C$ independent of $x,\,y,\,i$ such that
$g\equiv \wz C r_i^{1-s}\chi_{B(x_i,\,2C_0d(x_i,\,\boz^\complement))}$ is an element of $\cd^{s,\,1/(8C_0)}_\ball(u)$,
where $r_i=\diam S_i\sim d(x_i,\,\boz^\complement) $ by Definition \ref{d2.4}(iii), and
which implies that $u\in\dot M^{s,\,p}_\ball(\boz)$ and
$\|u_i\|_{\dot M^{s,\,p}_\ball(\boz)}\ls r_i^{1-s+n/p}$.
Notice also that $|u_i(x)-u_i(y)|\ge \dz_i$.

Moreover, let
\begin{equation}\label{e4.12}
 u\equiv \sum_{i=1}^{j-1}r_i^{(s-n/p)/(p-1)}r_i^{s-1-n/p}u_i
\end{equation}
and
$$g\equiv  \sum_{i=1}^{j-1} r_i^{(s-n/p)/(p-1)}r_i^{s-1-n/p}g_i= \sum_{i=1}^{j-1} r_i^{(s-n/p)/(p-1)}r_i^{-n/p}
\chi_{B(x_i,\,2C_0d(x_i,\,\boz^\complement))}.$$
Then there exists a positive constant $\wz C$ independent of $x,\,y$ such that $\overline Cg$ is
an element of $\cd^{s,\,1/(8C_0)}_\ball(u)$, which together with Definition \ref{d2.4} and
the vector-valued inequality of the
Hardy-Littlewood maximal function $\cm$ (see, for example, \cite{s93}) implies that $u\in\dot M^{s,\,p}_\ball(\boz)$ and
\begin{eqnarray*}
 \|u\|^{p}_{\dot M^{s,\,p}_\ball(\boz)}&&\ls
\int_\boz\lf(\sum_{i=1}^{j-1}r_i^{(s-n/p)/(p-1)}r_i^{-n/p}
\chi_{B(x_i,\,2C_0d(x_i,\,\boz^\complement))}(z)\r)^{p}\,dz\\
&&\ls
 \int_\boz\lf\{\sum_{i=1}^{j-1}
\lf[\cm \lf(\lf[r_i^{(s-n/p)/(p-1)}r_i^{-n/p}\chi_{B(x_i,\,C_0^{-1}d(x_i,\,\boz^\complement))}\r]^{1/2}\r)(z)\r]^{2}\r\}^{p}\,dz\nonumber\\
&&\ls
\int_\boz\lf(\sum_{i=1}^{j-1}r_i^{(s-n/p)/(p-1)}r_i^{-n/p}
 \chi_{B(x_i,\,C_0^{-1}d(x_i,\,\boz^\complement))} (z)\r)^{p}\,dz\nonumber\\
&&\ls \sum_{i=1}^ {j-1}r_i^{p(s-n/p)/(p-1)}.\nonumber
\end{eqnarray*}
By this and the assumption that $\boz$ is an $\dot M^{s,\,p}_\ball$-imbedding domain,
we have
\begin{eqnarray*}|u(x)-u(y)|
&&\ls |x-y|^{s-n/p}\|u\|_{ \dot M^{s,\,p}_\ball(\boz)}
\ls |x-y|^{s-n/p}\lf(\sum_{i=1}^{j-1}r_i^{p(s-n/p)/(p-1)}\r)^{1/p},
\end{eqnarray*}
which together with
$$\sum_{i=1}^{j-1}r_i^{(s-n/p)/(p-1)}r_i^{s-n/p} \le
\sum_{i=1}^{j-1}r_i^{(s-n/p)/(p-1)}r_i^{s-1-n/p}|u_i(x)-u_i(y)|\le |u(x)-u(y)|$$
implies that
$$ \sum_{i=1}^{j-1}r_i^{(ps-n)/(p-1)} \ls|x-y|^{(ps-n)/(p-1)}.$$
Thus
$$\int_\gz [d(z,\,\boz^\complement)]^{(ps-n)/(p-1)-1}\,|dz|\ls|x-y|^{(ps-n)/(p-1)},$$
which gives Theorem \ref{t1.2}(ii). This finishes the proof of Theorem \ref{t1.2}.
\end{proof}

\begin{rem}\label{r4.3}\rm
 Observe that in the proof of Theorem \ref{t1.2}(ii),
the functions $\{u_i\}_{i=1}^{j-1}$ and
hence $u$ defined in \eqref{e4.12} belong to $\dot M^{s,\,\fz}_\ball(\boz)$.
So to obtain that $\boz$ is a weak $(ps-n)/(p-1)$-cigar domain,
the assumption that $\boz$ is an $\dot M^{s,\,p}_\ball$-imbedding domain required in Theorem \ref{t1.2}(ii)
can be reduced to a weaker one: for all $u\in \dot M^{s,\,\fz}_\ball(\boz)$,
$\|u\|_{\dot M^{s-n/p,\,\fz}(\boz)}\ls\|u\|_{\dot M^{s ,\,p}_\ball(\boz)}$.
This leads to the proof of Theorem \ref{t1.3}.
\end{rem}

\begin{proof}[Proof of Theorem \ref{t1.3}.]
Let $\boz$ be a weak $\az$-cigar domain satisfying the slice property.
Let $\az^\ast\in(0,\,\az)$, $s\in(0,\,1]$ and $p^\ast,\,\wz p\in(n/s,\,\fz)$
be as in the proof of Theorem \ref{t1.2}(i).
Then by the conclusion \eqref{e4.9} there,
$\|u\|_{\dot M^{s-n/ p^\ast,\,\fz}(\boz)}\ls\|u\|_{\dot M^{s,\,  p^\ast}_\ball(\boz)}$
for all $u\in  M^{s,\, \wz p}_\ball(\boz)$ and hence for all $u\in \dot M^{s,\,\fz}_\ball(\boz)$.
So keeping the the observation in Remark \ref{r4.3} in mind and running the proof of Theorem \ref{t1.2}(ii) again,
 we obtain that $\boz$ is a weak $\az^\ast$-cigar domain.
Then by Theorem \ref{t1.2}(i), we further know that
$\boz$ is an $\dot M^{t,\,q}_\ball$-extension/-imbedding domain for all $t\in(\az^\ast,\,1]$
and $q\in[(n-\az^\ast)/(t-\az^\ast),\,\fz)$,
which completes the proof of Theorem \ref{t1.3}.
\end{proof}

{\bf Acknowledgements.} The author would like to thank Professor Pekka Koskela
for his kind suggestions and significant discussions
on this topic and also thank Professor Dachun Yang for his kind suggestions.

\noindent Yuan Zhou

\noindent Department of Mathematics and Statistics,
P. O. Box 35 (MaD),
FI-40014, University of Jyv\"askyl\"a,
Finland

\smallskip

\noindent{\it E-mail address}:  \texttt{yuan.y.zhou@jyu.fi}
\end{document}